\documentclass{article}
\usepackage{graphicx}
\usepackage{amsmath}
\usepackage{amsfonts}
\usepackage{amssymb}
\newtheorem{theorem}{Theorem}[subsection]

\newtheorem{corollary}[theorem]{Corollary}

\newtheorem{definition}[theorem]{Definition}
\newtheorem{example}[theorem]{Example}

\newtheorem{lemma}[theorem]{Lemma}

\newtheorem{proposition}[theorem]{Proposition}

\newenvironment{proof}[1][Proof]{\textbf{#1.} }{\ \rule{0.5em}{0.5em}}

\begin{document}

\title{ Combinatorics of crystal graphs and Kostka-Foulkes polynomials for the root
systems $B_{n},C_{n}$ and $D_{n}.$}
\author{C\'{e}dric Lecouvey\\lecouvey@math.unicaen.fr}
\date{}
\maketitle
\begin{abstract}
We use Kashiwara-Nakashima's combinatorics of crystal graphs associated to the
roots sytems $B_{n}$ and $D_{n}$ to extend the results of \cite{lec3} and
\cite{Mor} by showing that Morris type recurrence formulas also exist for the
orthogonal root systems.\ We derive from these formulas\ a statistic on
Kashiwara-Nakashima's tableaux of types $B_{n},C_{n}$ and $D_{n}$ generalizing
Lascoux-Sch\"{u}tzenberger's charge and from which it is possible to compute
the Kostka-Foulkes polynomials $K_{\lambda,\mu}(q)$ with restrictive
conditions on $(\lambda,\mu)$ . This statistic is different from that obtained
in \cite{lec3} from the cyclage graph structure on tableaux of type $C_{n}$.
We show that such a structure also exists for the tableaux of types $B_{n}$
and $D_{n}$ but can not be simply related to the Kostka-Foulkes polynomials.
Finally we give explicit formulas for $K_{\lambda,\mu}(q)$ when $\left|
\lambda\right|  \leq3,$ or $n=2$ and $\mu=0$.
\end{abstract}

\section{Introduction}

The multiplicity $K_{\lambda,\mu}$ of the weight $\mu$ in the irreducible
finite dimensional representation $V(\lambda)$ of the simple Lie algebra $g$
can be written in terms of the ordinary Kostant's partition function
$\mathcal{P}$ defined from the equality:%
\[
\prod_{\alpha\text{ positive root}}\dfrac{1}{(1-x^{\alpha})}=\sum_{\beta
}\mathcal{P}(\beta)x^{\beta}%
\]
where $\beta$ runs on the set of nonnegative integral combinations of positive
roots of $g$. Thus $\mathcal{P}(\beta)$ is the number of ways the weight
$\beta$ can be expressed as a sum of positive roots. Then we have%
\[
K_{\lambda,\mu}=\sum_{\sigma\in W}(-1)^{l(\sigma)}\mathcal{P}(\sigma
(\lambda+\rho)-(\mu+\rho))
\]
where $W$ is the Weyl group of $g.$

\noindent There exists a $q$-analogue $K_{\lambda,\mu}(q)$ of $K_{\lambda,\mu
}$ obtained by substituting the ordinary Kostant's partition function
$\mathcal{P}$ by its $q$-analogue $\mathcal{P}_{q}$ satisfying%
\[
\prod_{\alpha\text{ positive root}}\dfrac{1}{(1-qx^{\alpha})}=\sum_{\beta
}\mathcal{P}_{q}(\beta)x^{\beta}.
\]
So we have%
\[
K_{\lambda,\mu}(q)=\sum_{\sigma\in W}(-1)^{l(\sigma)}\mathcal{P}_{q}%
(\sigma(\lambda+\rho)-(\mu+\rho)).
\]
As shown by Lusztig \cite{Lut} $K_{\lambda,\mu}(q)$ is a polynomial in $q$
with non negative integer coefficients.

\noindent For type $A_{n-1}$ the positivity of the Kostka-Foulkes Polynomials
can also be proved by a purely combinatorial method. Recall that for any
partitions $\lambda$ and $\mu$ with $n$ parts the number of semi-standard
tableaux of shape $\lambda$ and weight $\mu$ is equal to the multiplicity of
the weight $\mu$ in $V(\lambda).$

In \cite{LSc1}, Lascoux and Sch\"{u}tzenberger have introduced a beautiful
statistic $\mathrm{ch}_{A}$ on dominant evaluation words $w$ that is on words
$w=x_{1}\cdot\cdot\cdot x_{l}$ whose letters $x_{i}$ are positive integers
such that for any $i\geq1$ with $i$ a letter of $w,$ $w$ contains more letters
$i$ than letters $i+1.$ Recall that the plactic monoid is the quotient set of
the free monoid on the positive integers by Knuth's relations%
\[
abx\equiv\left\{
\begin{tabular}
[c]{l}%
$bax$ if $a<x\leq b$\\
$axb$ if $x\leq a<b$%
\end{tabular}
\right.  .
\]
The statistic $\mathrm{ch}_{A}$ is the unique function from dominant
evaluation words to non-negative integers such that
\begin{equation}
\left\{
\begin{tabular}
[c]{l}%
$\mathrm{ch}_{A}(\emptyset)=0$\\
$\mathrm{ch}_{A}(xu)=\mathrm{ch}_{A}(ux)+1$ if $x$ is not the lowest letter of
$w$\\
$\mathrm{ch}_{A}(xu)=\mathrm{ch}_{A}(u)$ if $x$ is the lowest letter of $w$\\
$\mathrm{ch}_{A}(\sigma w)=\mathrm{ch}_{A}(w)$ for any $\sigma\in
\mathcal{S}_{n}$\\
$\mathrm{ch}_{A}(w_{1})=\mathrm{ch}_{A}(w_{2})$ if $w_{1}\equiv w_{2}$%
\end{tabular}
\right.  \label{def_chA_ini}%
\end{equation}
\cite{kill}. Then the charge of the semi-standard tableau $T$ of dominant
weight verifies $\mathrm{ch}_{A}(T)=\mathrm{ch}_{A}(\mathrm{w}(T))$ where
$\mathrm{w}(T)$ is the word obtained by column reading the letters of $T$ from
top to bottom and right to left. Lascoux and Sch\"{u}tzenberger have proved
the equality%
\begin{equation}
K_{\lambda,\mu}(q)=\sum_{T}q^{\mathrm{ch}_{A}(T)} \label{th_ls}%
\end{equation}
where $T$ runs on the set of semi-standard tableaux of shape $\lambda$ and
weight $\mu.\;$The proof of (\ref{th_ls}) is based on Morris recurrence
formula which permits to express each Kostka-Foulkes polynomials related to
the root system $A_{n}$ in terms of Kostka-Foulkes polynomials related to the
root system $A_{n-1}.$

\noindent The compatibility of the charge with plactic relations provides
alternative ways to compute $\mathrm{ch}_{A}(T)$. By applying the reverse
bumping algorithm on the boxes contained in the longest row of $T$ we obtain a
pair $(R,T^{\prime})$ with $R$ a row tableau whose length is equal to the
longest row of $T$ and $T^{\prime}$ a semi-standard tableau which does not
contain the lowest letter $t$ of $T$ such that $\mathrm{w}(T)\equiv
\mathrm{w}(R)\otimes\mathrm{w}(T^{\prime})$. Let $R^{\prime}$ be the row
tableau obtained by erasing all the letters $t$ in $R.$ Then the catabolism of
$T$ is the unique semi-standard tableau $\mathrm{cat}(T)$ such that
$\mathrm{w}(\mathrm{cat}(T))\equiv\mathrm{w}(T^{\prime})\otimes\mathrm{w}%
(R^{\prime})$ computed via the bumping algorithm. We have
\[
\mathrm{ch}_{A}(\mathrm{cat}(T))=\mathrm{ch}_{A}(T)+r^{\prime}%
\]
where $r^{\prime}$ is the length of $R^{\prime}.$ Since the number of boxes of
$\mathrm{cat}(T)$ is strictly less than that of $T,$ $\mathrm{ch}_{A}(T)$ can
be obtained from $T$ by computing successive catabolism operations$.$ In fact
this is this characterization of the charge which is needed to prove
(\ref{th_ls}).

\noindent The charge may also be obtained by endowing $ST(\mu)$ the set of
semi-standard tableaux of weight $\mu$ with a structure of graph. We draw an
arrow $T\rightarrow T^{\prime}$ between the two tableaux $T$ and $T^{\prime}$
of $ST(\mu),$ if and only if there exists a word $u$ and a letter $y$ which is
not the lowest letter of $T$ such that $\mathrm{w}(T)\equiv xu$ and
$\mathrm{w}(T^{\prime})\equiv ux$.\ Then we say that $T^{\prime}$ is a
cocyclage of $T.$ The essential tool to define this graph structure is yet the
bumping algorithm for the semi-standard tableaux. The cyclage graph $ST(\mu)$
contains a unique row tableau $L_{\mu}$ which can not be obtained as the
cocyclage of another tableau of $ST(\mu)$. Let $T_{\mu}$ be the unique
semi-standard tableau of shape $\mu$ belonging to $ST(\mu)$. Then there is no
cocyclage of $T_{\mu}.$ For any $T\in ST(\mu)$ all the paths joining $L_{\mu}$
to $T$ have the same length.\ This length is called the cocharge of $T$ and
denoted $\mathrm{coch}_{A}(T).$ Similarly, all the paths joining $T$ to
$T_{\mu}$ have the same length which is equal to the charge of $T.$ The
maximal value of $\mathrm{ch}_{A}$ is $\left\|  \mu\right\|  =\mathrm{ch}%
_{A}(L_{\lambda})=\sum_{i}(i-1)\mu_{i}.\;$Moreover the charge and the cocharge
are related by the equality $\mathrm{ch}_{A}(T)=\left\|  \mu\right\|
-\mathrm{coch}_{A}(T)$ for any $T\in ST(\mu).$

\bigskip

Analogues of semi-standard tableaux also exit for the other classical root
systems. They have been introduced by Kashiwara and Nakashima \cite{KN} via
crystal bases theory. For each classical root system these tableaux naturally
label the vertices of the crystal graph $B(\lambda)$ associated to the
dominant weight $\lambda.$

\noindent In \cite{lec3} we have proved that an analogue of Morris recurrence
formula exists for the root system $C_{n}.$ Moreover it is also possible to
endow the corresponding set of tableaux with a structure of cyclage graph.
From these graphs we have introduced a natural statistic on
Kashiwara-Nakashima's tableaux of type $C_{n}$ and have conjectured that this
statistic yields an analogue of Lascoux-Sch\"{u}tzenberger's theorem.

\bigskip

This article is an attempt to look at possible generalizations and extensions
of these results to the orthogonal roots systems. We establish Morris type
recurrence formula for the root systems $B_{n}$ and $D_{n}.$ Moreover we show
that is possible to endow the set of tableaux of types $B_{n}$ and $D_{n}$
with a structure of cyclage graph. Nevertheless the situation is more
complicated than for the root system $C_{n}$ and we are not able to deduce
from these graphs a natural statistic relevant for computing the
Kostka-Foulkes polynomials. To overcome this problem we change our strategy
and define a new statistic $\chi_{n}$ on tableaux of types $B_{n},C_{n}$ and
$D_{n}$ based on the catabolism operation. Then we prove that this statistic
can be used to compute the Kostka-Foulkes polynomials $K_{\lambda,\mu}(q)$
with restrictive conditions on $(\lambda,\mu)$. Note that the analogue of
(\ref{th_ls}) with $\chi_{n}$ is false in general. In particular $\chi_{n}$ is
not equal to the statistic defined in \cite{lec3} for the tableaux of type
$C_{n}$ even if the two statistics can be regarded as generalizations of
$\mathrm{ch}_{A}$ since they coincide on semi-standard tableaux.

\bigskip

In Section $1$ we recall the Background on Kostka-Foulkes polynomials and
combinatorics of crystal graphs that we need in the sequel.\ We also summarize
the basic properties of the insertion algorithms and plactic monoids for the
root systems $B_{n},C_{n}$ and $D_{n}$ introduced in \cite{Lec} and
\cite{lec2}.\ Section $2$ is devoted to Morris type recurrence formulas for
types $B_{n}$ and $C_{n}$. In Section $3$ we define the catabolism operation
for the tableaux of type $B_{n},C_{n}$ and $D_{n}.$ Then we introduce the
statistics $\chi_{n}^{B},\chi_{n}^{C}$ and $\chi_{n}^{D}$ and prove that
analogues of (\ref{th_ls}) hold for these statistics if $\lambda$ and $\mu$
satisfy restrictive conditions. We also introduce the cyclage graph structure
on tableaux of types $B_{n}$ and $D_{n}$ and show that a charge statistic
related to Kostka-Foulkes polynomials can not be obtained in a similar way
that in \cite{lec3}. Finally we give in Section $4$ explicit simple formulas
for the Kostka-Foulkes polynomials $K_{\lambda,\mu}(q)$ when $\left|
\lambda\right|  \leq3,$ or $n=2$ and $\mu=0$ deduced from the results of
Sections $2$ and $3.$

\bigskip

\noindent\textbf{Notation: }In the sequel we frequently define similar objects
for the root systems $B_{n}$ $C_{n}$ and $D_{n}$. When they are related to
type $B_{n}$ (resp. $C_{n},D_{n}$), we implicitly attach to them the label $B$
(resp. the labels $C,D$). To avoid cumbersome repetitions, we sometimes omit
the labels $B,C$ and $D$ when our definitions or statements are identical for
the three root systems.

\bigskip

\section{Background}

\subsection{Kostka-Foulkes polynomials associated to a root system}

Let $g$ be a simple Lie algebra and $\alpha_{i},$ $i\in I$ its simple roots.
Write $Q^{+}$ and $R^{+}$ for the set of nonnegative integral combinations of
positive roots and for the set of positive roots of $g$. Denote respectively
by $P$ and $P^{+}$ its weight lattice and its cone of dominant weights. Let
$\{s_{i},i\in I\}$ be a set of generators of the Weyl group $W$ and $l$ the
corresponding length function.

\noindent The $q$-analogue $\mathcal{P}_{q}$ of the Kostant function partition
is such that%
\[
\prod_{\alpha\in R^{+}}\dfrac{1}{1-qx^{\alpha}}=\sum_{\beta\in Q^{+}%
}\mathcal{P}_{q}(\beta)x^{\beta}\text{ and }\mathcal{P}_{q}(\beta)=0\text{ if
}\beta\notin Q^{+}.
\]

\begin{definition}
\label{def_kost}Let $\lambda,\mu\in P^{+}.$ The Kostka-Foulkes polynomial
$K_{\lambda,\mu}(q)$ is defined by
\[
K_{\lambda,\mu}(q)=\sum_{\sigma\in W}(-1)^{l(\sigma)}\mathcal{P}_{q}%
(\sigma(\lambda+\rho)-(\mu+\rho)).
\]
where $\rho$ is the half sum of positive roots.
\end{definition}

Let $\beta\in P.\;$We set
\[
a_{\beta}=\sum_{\sigma\in W}(-1)^{l(\sigma)}(\sigma\cdot x^{\beta})
\]
where $\sigma\cdot x^{\mu}=x^{\sigma(\mu)}.$ The Schur function $s_{\beta}$ is
defined by
\[
s_{\beta}=\dfrac{a_{\beta+\rho}}{a_{\rho}}.
\]
When $\lambda\in P^{+},$ $s_{\lambda}$ is the Weyl character of $V(\lambda)$
the finite dimensional irreducible $g$-module with highest weight $\lambda.$
For any $\sigma\in W,$ the dot action of $\sigma$ on $\beta\in P$ is defined
by $\sigma\circ\beta=\sigma\cdot(\beta+\rho)-\rho.$ We have the following
straightening law for the Schur functions. For any $\beta\in P$, $s_{\beta}=0$
or there exists a unique $\lambda\in P^{+}$ such that $s_{\beta}%
=(-1)^{l(\sigma)}s_{\lambda}$ with $\sigma\in W$ and $\lambda=\sigma\circ
\beta.$ Set $\mathbb{K}=\mathbb{Z}[q,q^{-1}]$ and write $\mathbb{K}[P]$ for
the $\mathbb{K}$-module generated by the $x^{\beta}$, $\beta\in P.$ Set
$\mathbb{K}[P]^{W}=\{f\in\mathbb{K}[P],$ $\sigma\cdot f=f$ for any $\sigma\in
W\}.$ Then $\{s_{\lambda}\}$ is a basis of $\mathbb{K}[P]^{W}.$

\noindent To each positive root $\alpha,$ we associate the raising operator
$R_{\alpha}:P\rightarrow P$ defined by%
\[
R_{\alpha}(\beta)=\alpha+\beta.
\]
Given $\alpha_{1},...,\alpha_{p}$ positive roots and $\beta\in P,$ we set
$(R_{\alpha_{1}}\cdot\cdot\cdot R_{\alpha_{p}})s_{\beta}=s_{R_{\alpha_{1}%
}\cdot\cdot\cdot R_{\alpha_{p}}(\beta)}.$ For all $\beta\in P,$ we define the
Hall-Littelwood polynomial $Q_{\beta}$ by%
\[
Q_{\beta}=\left(  \prod_{\alpha\in R^{+}}\dfrac{1}{1-qR_{\alpha}}\right)
s_{\beta}%
\]
where $\dfrac{1}{1-qR_{\alpha}}=\sum_{k=0}^{+\infty}q^{k}R_{\alpha}^{k}.$

\begin{theorem}
\label{th_hall_kostka}\cite{mac}For any $\lambda,\mu\in P^{+},$ $K_{\lambda
,\mu}(q)$ is the coefficient of $s_{\lambda}$ in $Q_{\mu}$ that is,
\[
Q_{\mu}=\sum_{\lambda\in P^{+}}K_{\lambda,\mu}(q)s_{\lambda}.
\]
\end{theorem}

\subsection{Kostka-Foulkes polynomials for the root systems $B_{n},C_{n}$ and
$D_{n}\label{subsec_KF}$}

We choose to label respectively the Dynkin diagrams of $so_{2n+1}$, $sp_{2n}$
and $so_{2n}$ by
\begin{equation}
\overset{0}{\circ}\Longleftarrow\overset{1}{\circ}-\overset{2}{\circ}%
-\overset{3}{\circ}-\overset{4}{\circ}-\cdot\cdot\cdot\overset{n-1}{\circ
}\text{, }\overset{0}{\circ}\Longrightarrow\overset{1}{\circ}-\overset
{2}{\circ}-\overset{3}{\circ}-\overset{4}{\circ}-\cdot\cdot\cdot\overset
{n-1}{\circ}\text{ and }%
\begin{tabular}
[c]{l}%
$\overset{1}{\circ}$\\
$\ \ \backslash$\\
$\ \ \ \ \overset{2}{\circ}$\\
$\ \ /$\\
$\underset{0}{\circ}$%
\end{tabular}
-\overset{3}{\circ}-\overset{4}{\circ}-\cdot\cdot\cdot\overset{n-2}{\circ
}-\overset{n-1}{\circ}. \label{DD}%
\end{equation}
The weight lattices for the root systems $B_{n},C_{n}$ and $D_{n}$ can be
identified with $P_{n}=\mathbb{Z}^{n}$ equipped with the orthonormal basis
$\varepsilon_{\overline{i}},$ $i=1,...,n$.\ We take for the simple roots%
\begin{equation}
\left\{
\begin{tabular}
[c]{l}%
$\alpha_{0}^{B_{n}}=\varepsilon_{\overline{1}}\text{ and }\alpha_{i}^{B_{n}%
}=\varepsilon_{\overline{i+1}}-\varepsilon_{\overline{i}}\text{,
}i=1,...,n-1\text{ for the root system }B_{n}$\\
$\alpha_{0}^{C_{n}}=2\varepsilon_{\overline{1}}\text{ and }\alpha_{i}^{C_{n}%
}=\varepsilon_{\overline{i+1}}-\varepsilon_{\overline{i}}\text{,
}i=1,...,n-1\text{ for the root system }C_{n}$\\
$\alpha_{0}^{D_{n}}=\varepsilon_{\overline{1}}+\varepsilon_{\overline{2}%
}\text{ and }\alpha_{i}^{D_{n}}=\varepsilon_{\overline{i+1}}-\varepsilon
_{\overline{i}}\text{, }i=1,...,n-1\text{ for the root system }D_{n}$%
\end{tabular}
\right.  . \label{simple_roots}%
\end{equation}
Then the set of positive roots are%
\[
\left\{
\begin{tabular}
[c]{l}%
$R_{B_{n}}^{+}=\{\varepsilon_{\overline{i}}-\varepsilon_{\overline{j}%
},\varepsilon_{\overline{i}}+\varepsilon_{\overline{j}}\text{ with }1\leq
j<i\leq n\}\cup\{\varepsilon\overline{_{i}}\text{ with }1\leq i\leq n\}\text{
for the root system }B_{n}$\\
$R_{B_{n}}^{+}=\{\varepsilon_{\overline{i}}-\varepsilon_{\overline{j}%
},\varepsilon_{\overline{i}}+\varepsilon_{\overline{j}}\text{ with }1\leq
j<i\leq n\}\cup\{2\varepsilon\overline{_{i}}\text{ with }1\leq i\leq n\}\text{
for the root system }C_{n}$\\
$R_{D_{n}}^{+}=\{\varepsilon_{\overline{i}}-\varepsilon_{\overline{j}%
},\varepsilon_{\overline{i}}+\varepsilon_{\overline{j}}\text{ with }1\leq
j<i\leq n\}\text{ for the root system \ }D_{n}$%
\end{tabular}
\right.  .
\]
Denote respectively by $P_{B_{n}}^{+},P_{C_{n}}^{+}$ and $P_{D_{n}}^{+}$the
sets of dominant weights of $so_{2n+1},sp_{2n}$ and $so_{2n}.$ Write
$\Lambda_{0}^{B_{n}},...,\Lambda_{n-1}^{B_{n}}$ for the fundamentals weights
of $so_{2n+1},$ $\Lambda_{0}^{C_{n}},...,\Lambda_{n-1}^{C_{n}}$ for the
fundamentals weights of $sp_{2n}$ and $\Lambda_{0}^{D_{n}},...,\Lambda
_{n-1}^{D_{n}}$ for the fundamentals weights of $so_{2n+1}$.

\noindent We have $\Lambda_{i}^{B_{n}}=\Lambda_{i}^{C_{n}}=\Lambda_{i}^{D_{n}%
}=\varepsilon_{\overline{n}}+\cdot\cdot\cdot+\varepsilon_{\overline{i+1}}$ for
$2\leq i\leq n-1$, $\Lambda_{0}^{B_{n}}=\Lambda_{0}^{D_{n}}=\dfrac{1}%
{2}(\varepsilon_{\overline{n}}+\cdot\cdot\cdot+\varepsilon_{\overline{2}%
}+\varepsilon_{\overline{1}})$, $\Lambda_{0}^{C_{n}}=\varepsilon_{\overline
{n}}+\cdot\cdot\cdot+\varepsilon_{\overline{2}}+\varepsilon_{\overline{1}},$
$\Lambda_{1}^{B_{n}}=\Lambda_{1}^{C_{n}}=\varepsilon_{\overline{n}}+\cdot
\cdot\cdot+\varepsilon_{\overline{2}}$ and $\Lambda_{1}^{D_{n}}=\dfrac{1}%
{2}(\varepsilon_{\overline{n}}+\cdot\cdot\cdot+\varepsilon_{\overline{2}%
}-\varepsilon_{\overline{1}}).$

\noindent Consider $\lambda\in P_{B_{n}}^{+}$ and write $\lambda=\sum
_{i=0}^{n-1}\widehat{\lambda}_{i}\Lambda_{i}^{B}$ with $\widehat{\lambda}%
_{i}\in\mathbb{N}.\;$Set $\lambda_{\overline{1}}=\dfrac{\widehat{\lambda}_{0}%
}{2}$ and $\lambda_{\overline{i}}=\dfrac{\widehat{\lambda}_{0}}{2}%
+\widehat{\lambda}_{1}+\cdot\cdot\cdot+\widehat{\lambda}_{i-1},$ $i=2,...,n.$
The dominant weight $\lambda$ is characterized by the generalized partition
$(\lambda_{\overline{n}},...,\lambda_{\overline{1}})$ such that $\lambda
_{\overline{n}}\geq\cdot\cdot\cdot\geq\lambda_{\overline{1}}$ and
$\lambda_{\overline{i}}\in\dfrac{\mathbb{N}}{2},$ $i=1,...,n.$ In the sequel
we will identify $\lambda$ and $(\lambda_{\overline{n}},...,\lambda
_{\overline{1}})$ by setting $\lambda=(\lambda_{\overline{n}},...,\lambda
_{\overline{1}}).$ Then $\lambda=\lambda_{\overline{1}}\varepsilon
_{\overline{1}}+\cdot\cdot\cdot+\lambda_{\overline{n}}\varepsilon
_{\overline{n}}$ that is, the $\lambda_{i}$'s are the coordinates of $\lambda$
on the basis $(\varepsilon_{\overline{n}},...,\varepsilon_{\overline{1}}%
)$.\ The half sum of positive roots verifies $\rho_{B_{n}}=(n-\dfrac{1}%
{2},n-\dfrac{3}{2},...,\dfrac{1}{2}).$

\noindent Consider $\lambda\in P_{C_{n}}^{+}$ and write $\lambda=\sum
_{i=0}^{n-1}\widehat{\lambda}_{i}\Lambda_{i}^{C}$ with $\widehat{\lambda}%
_{i}\in\mathbb{N}.\;$The dominant weight $\lambda$ is characterized by the
partition $(\lambda_{\overline{n}},...,\lambda_{\overline{1}})$ where
$\lambda_{\overline{1}}=\widehat{\lambda}_{0}$ and $\lambda_{\overline{i}%
}=\widehat{\lambda}_{0}+\widehat{\lambda}_{1}+\cdot\cdot\cdot+\widehat
{\lambda}_{i-1},$ $i=2,...,n.$ We set $\lambda=(\lambda_{\overline{n}%
},...,\lambda_{\overline{1}}).$ Then $\lambda=\lambda_{\overline{1}%
}\varepsilon_{\overline{1}}+\cdot\cdot\cdot+\lambda_{\overline{n}}%
\varepsilon_{\overline{n}}$ and the half sum of positive roots verifies
$\rho_{C_{n}}=(n,n-1,...,1).$

\noindent Now consider $\lambda\in P_{D_{n}}^{+}$ and write $\lambda
=\sum_{i=0}^{n-1}\widehat{\lambda}_{i}\Lambda_{i}^{D}$ with $\widehat{\lambda
}_{i}\in\mathbb{N}.\;$Set $\lambda_{\overline{1}}=\dfrac{\widehat{\lambda}%
_{0}-\widehat{\lambda}_{1}}{2}$, $\lambda_{\overline{2}}=\dfrac{\widehat
{\lambda}_{0}+\widehat{\lambda}_{1}}{2}$ and $\lambda_{\overline{i}}%
=\dfrac{\widehat{\lambda}_{0}+\widehat{\lambda}_{1}}{2}+\widehat{\lambda}%
_{2}+\cdot\cdot\cdot+\widehat{\lambda}_{i-1},$ $i=3,...,n.$ The dominant
weight $\lambda$ is characterized by the generalized partition $(\lambda
_{\overline{n}},...,\lambda_{\overline{1}})$ such that $\lambda_{\overline{n}%
}\geq\cdot\cdot\cdot\geq\lambda_{\overline{1}}$, $\lambda_{\overline{i}}%
\in\dfrac{\mathbb{N}}{2}$ $i=2,...,n$ and $\lambda_{\overline{1}}\in
\dfrac{\mathbb{Z}}{2}.$ Note that we can have $\lambda_{\overline{1}}<0.$ We
set $\lambda=(\lambda_{\overline{n}},...,\lambda_{\overline{1}}).$ Then
$\lambda=\lambda_{\overline{1}}\varepsilon_{\overline{1}}+\cdot\cdot
\cdot+\lambda_{\overline{n}}\varepsilon_{\overline{n}}$ and the half sum of
positive roots verifies $\rho_{D_{n}}=(n-1,n-2,...,0).$

\bigskip

For any generalized partition $\lambda=(\lambda_{\overline{n}},...,\lambda
_{\overline{1}})\in P_{n}^{+},$ we write $\lambda^{\prime}\in P_{n-1}^{+}$ for
the generalized partition obtained by deleting $\lambda_{\overline{n}}$ in
$\lambda.$ Moreover we set $\left|  \lambda\right|  =\lambda_{\overline{1}%
}+\lambda_{\overline{2}}+\cdot\cdot\cdot+\lambda_{\overline{n}}$ if
$\lambda_{\overline{1}}\geq0,$ $\left|  \lambda\right|  =-\lambda
_{\overline{1}}+\lambda_{\overline{2}}+\cdot\cdot\cdot+\lambda_{\overline{n}}$ otherwise.

\noindent The Weyl group $W_{B_{n}}=W_{C_{n}}$ of $so_{2n+1}$ can be regarded
as the sub group of the permutation group of $\{\overline{n},...,\overline
{2},\overline{1},1,2,...,n\}$\ generated by $s_{i}=(i,i+1)(\overline
{i},\overline{i+1}),$ $i=1,...,n-1$ and $s_{0}=(1,\overline{1})$ where for
$a\neq b$ $(a,b)$ is the simple transposition which switches $a$ and $b.$ We
denote by $l_{B}$ the length function corresponding to the set of generators
$s_{i},$ $i=0,...n-1.$

\noindent The Weyl group $W_{D_{n}}$ of $so_{2n}$ can be regarded as the sub
group of the permutation group of $\{\overline{n},...,\overline{2}%
,\overline{1},1,2,...,n\}$\ generated by $s_{i}=(i,i+1)(\overline{i}%
,\overline{i+1}),$ $i=1,...,n-1$ and $s_{0}^{\prime}=(1,\overline
{2})(2,\overline{1})$. We denote by $l_{D}$ the length function corresponding
to the set of generators $s_{0}^{\prime}$ and $s_{i},$ $i=1,...n-1.$

\noindent Note that $W_{D_{n}}\subset W_{B_{n}}$ and any $\sigma\in W_{B_{n}}$
verifies $\sigma(\overline{i})=\overline{\sigma(i)}$ for $i\in\{1,...,n\}.$
The action of $\sigma$ on $\beta=(\beta_{\overline{n}},...,\beta_{\overline
{1}})\in P_{n}$ is given by%
\[
\sigma\cdot(\beta_{\overline{n}},...,\beta_{\overline{1}})=(\beta
_{\overline{n}}^{\sigma},...,\beta_{\overline{1}}^{\sigma})
\]
where $\beta_{\overline{i}}^{\sigma}=\beta_{\sigma(\overline{i})}$ if
$\sigma(\overline{i})\in\{\overline{1},...,\overline{n}\}$ and $\beta
_{\overline{i}}^{\sigma}=-\beta_{\sigma(i)}$ otherwise.

\noindent For any $\beta=(\beta_{\overline{n}},...,\beta_{\overline{1}})\in
P_{n}$ we set $x^{\beta}=x_{n}^{\beta_{\overline{n}}}\cdot\cdot\cdot
x_{1}^{\beta_{\overline{1}}}$ where $x_{1},...,x_{n}$ are fixed indeterminates.

\noindent The following lemma is a consequence of Definition \ref{def_kost}.

\begin{proposition}
\label{prop_degreeK}The Kostka-Foulkes polynomial $K_{\lambda},_{\mu}(q)$ is
monic of degree

\begin{itemize}
\item $\sum_{i=1}^{n}i(\lambda\overline{_{i}}-\mu\overline{_{i}})$ for the
root system $B_{n}$

\item $\sum_{i=1}^{n}i(\lambda\overline{_{i}}-\mu\overline{_{i}})-\dfrac{1}%
{2}(\left|  \lambda\right|  -\left|  \mu\right|  )$ for the root system $C_{n}$

\item $\sum_{i=2}^{n}(i-1)(\lambda\overline{_{i}}-\mu\overline{_{i}})$ for the
root system $D_{n}$
\end{itemize}
\end{proposition}

\begin{proof}
It is similar to that given in Example 4 page 243 of \cite{mac} for the degree
of Kostka-Foulkes polynomials associated to the root system $A_{n}.$
\end{proof}

\bigskip

\noindent\textbf{Remarks:}

\noindent$\mathrm{(i):}$ The above proposition suffices to determinate
$K_{\lambda,\mu}(q)$ when $\dim V(\lambda)_{\mu}=1.$ In particular we have
$K_{\lambda,\mu}(q)=1$ for each minuscule representation $V(\lambda).$

\noindent$\mathrm{(ii):}$ If $\left|  \lambda\right|  =\left|  \mu\right|  $
then $K_{\lambda,\mu}^{B_{n}}(q)=K_{\lambda,\mu}^{C_{n}}(q)=K_{\lambda,\mu
}^{D_{n}}(q)=K_{\lambda,\mu}^{A_{n-1}}(q).$

\noindent$\mathrm{(iii):}$ Suppose $\lambda,\mu\in P_{D_{n}}^{+}$. Set
$\lambda^{\ast}=(\lambda_{\overline{n}},...,\lambda_{\overline{2}}%
,-\lambda_{\overline{1}})$ and $\mu^{\ast}=(\mu_{\overline{n}},...,\mu
_{\overline{2}},-\mu_{\overline{1}})$ then
\begin{equation}
K_{\lambda,\mu}^{D_{n}}(q)=K_{\lambda^{\ast},\mu^{\ast}}^{D_{n}}(q)
\label{K=K*}%
\end{equation}
This is due to the symmetric role played by the simple roots $\alpha_{0}$ and
$\alpha\overline{_{1}}$ in the root system $D_{n}.$ Moreover when
$\lambda_{\overline{1}}=0,$ $\lambda=\lambda^{\ast}$ thus $K_{\lambda,\mu
}^{D_{n}}(q)=K_{\lambda^{\ast},\mu^{\ast}}^{D_{n}}(q)=K_{\lambda^{\ast},\mu
}^{D_{n}}(q)=K_{\lambda,\mu^{\ast}}^{D_{n}}(q).$

\noindent$\mathrm{(iv):}$ Consider $\lambda,\mu\in P_{n}^{+}$ such that
$\lambda_{\overline{n}}=\mu_{\overline{n}}.$ Then $K_{\lambda,\mu
}(q)=K_{\lambda^{\prime},\mu^{\prime}}(q).$

\subsection{Convention for crystal graphs}

In the sequel $g$ is any of the Lie algebras $so_{2n+1},sp_{2n}$ or $so_{2n}$.
The crystal graphs for the $U_{q}(g)$-modules are oriented colored graphs with
colors $i\in\{0,...,n-1\}$. An arrow $a\overset{i}{\rightarrow}b$ means that
$\widetilde{f}_{i}(a)=b$ and $\widetilde{e}_{i}(b)=a$ where $\widetilde{e}%
_{i}$ and $\widetilde{f}_{i}$ are the crystal graph operators (for a review of
crystal bases and crystal graphs see \cite{Ka2}). A vertex $v^{0}\in B$
satisfying $\widetilde{e}_{i}(v^{0})=0$ for any $i\in\{0,...,n-1\}$ is called
a highest weight vertex. The decomposition of $V$ into its irreducible
components is reflected into the decomposition of $B$ into its connected
components. Each connected component of $B$ contains a unique highest weight
vertex.\ The crystals graphs of two isomorphic irreducible components are
isomorphic as oriented colored graphs. The action of $\widetilde{e}_{i}$ and
$\widetilde{f}_{i}$ on $B\otimes B^{\prime}=\{b\otimes b^{\prime};$ $b\in
B,b^{\prime}\in B^{\prime}\}$ is given by:%

\begin{align}
\widetilde{f_{i}}(u\otimes v)  &  =\left\{
\begin{tabular}
[c]{c}%
$\widetilde{f}_{i}(u)\otimes v$ if $\varphi_{i}(u)>\varepsilon_{i}(v)$\\
$u\otimes\widetilde{f}_{i}(v)$ if $\varphi_{i}(u)\leq\varepsilon_{i}(v)$%
\end{tabular}
\right. \label{TENS1}\\
&  \text{and}\nonumber\\
\widetilde{e_{i}}(u\otimes v)  &  =\left\{
\begin{tabular}
[c]{c}%
$u\otimes\widetilde{e_{i}}(v)$ if $\varphi_{i}(u)<\varepsilon_{i}(v)$\\
$\widetilde{e_{i}}(u)\otimes v$ if $\varphi_{i}(u)\geq\varepsilon_{i}(v)$%
\end{tabular}
\right.  \label{TENS2}%
\end{align}
where $\varepsilon_{i}(u)=\max\{k;\widetilde{e}_{i}^{k}(u)\neq0\}$ and
$\varphi_{i}(u)=\max\{k;\widetilde{f}_{i}^{k}(u)\neq0\}$. The weight of the
vertex $u$ is defined by $\mathrm{wt}(u)=\underset{i=0}{\overset{n-1}{\sum}%
}(\varphi_{i}(u)-\varepsilon_{i}(u))\Lambda_{i}$.

\noindent The following lemma is a straightforward consequence of
(\ref{TENS1}) and (\ref{TENS2}).

\begin{lemma}
\label{lem_plu_hp}Let $u\otimes v$ $\in$ $B\otimes B^{\prime}$ $u\otimes v$ is
a highest weight vertex of $B\otimes B^{\prime}$ if and only if for any
$i\in\{0,...,n-1\}$ $\widetilde{e}_{i}(u)=0$ (i.e. $u$ is of highest weight)
and $\varepsilon_{i}(v)\leq\varphi_{i}(u).$
\end{lemma}

\noindent The Weyl group $W$ acts on $B$ by:
\begin{align}
s_{i}(u)  &  =(\widetilde{f_{i}})^{\varphi_{i}(u)-\varepsilon_{i}(u)}(u)\text{
if }\varphi_{i}(u)-\varepsilon_{i}(u)\geq0,\label{actionW}\\
s_{i}(u)  &  =(\widetilde{e_{i}})^{\varepsilon_{i}(u)-\varphi_{i}(u)}(u)\text{
if }\varphi_{i}(u)-\varepsilon_{i}(u)<0.\nonumber
\end{align}
We have the equality $\mathrm{wt}(\sigma(u))=\sigma(\mathrm{wt}(u))$ for any
$\sigma\in W$ and $u\in B.$ For any $\lambda\in P^{+},$ we denote by
$B(\lambda)$ the crystal graph of $V(\lambda).$

\subsection{Kashiwara-Nakashima's tableaux}

\noindent Accordingly to (\ref{DD}) the crystal graphs of the vector
representations are:%
\begin{gather*}
B(\Lambda_{n-1}^{B}):\overline{n}\overset{n-1}{\rightarrow}\overline
{n-1}\overset{n-2}{\rightarrow}\cdot\cdot\cdot\cdot\rightarrow\overline
{2}\overset{1}{\rightarrow}\overline{1}\overset{0}{\rightarrow}0\overset
{0}{\rightarrow}1\overset{1}{\rightarrow}2\cdot\cdot\cdot\cdot\overset
{n-2}{\rightarrow}n-1\overset{n-1}{\rightarrow}n\\
B(\Lambda_{n-1}^{C}):\overline{n}\overset{n-1}{\rightarrow}\overline
{n-1}\overset{n-2}{\rightarrow}\cdot\cdot\cdot\cdot\rightarrow\overline
{2}\overset{1}{\rightarrow}\overline{1}\overset{0}{\rightarrow}1\overset
{1}{\rightarrow}2\cdot\cdot\cdot\cdot\overset{n-2}{\rightarrow}n-1\overset
{n-1}{\rightarrow}n\\
B(\Lambda_{n-1}^{D}):\overline{n}\overset{n-1}{\rightarrow}\overline
{n-1}\overset{n-2}{\rightarrow}\cdot\cdot\cdot\overset{3}{\rightarrow
}\overline{3}\overset{2}{\rightarrow}%
\begin{tabular}
[c]{c}%
$1$ \ \ \\
\ \ $\overset{0}{\nearrow}$ $\ \ \ \overset{\text{ \ \ \ }1}{\text{
\ }\searrow}$ \ \ \ \\
$\overline{2}\ \ \ \ \ \ \ \ \ \ \ \ \ \ \ \ \ \ 2$\\
\ $\underset{1\text{ \ \ \ }}{\searrow}$ \ \ \ $\underset{0}{\nearrow}$
\ \ \ \\
$\overline{1}$ \
\end{tabular}
\overset{2}{\rightarrow}3\overset{3}{\rightarrow}\cdot\cdot\cdot\overset
{n-2}{\rightarrow}n-1\overset{n-1}{\rightarrow}n.
\end{gather*}
Kashiwara-Nakashima's combinatorial description of the crystal graphs
$B(\lambda)$ is based on a notion of tableaux analogous for each root system
$B_{n},C_{n}$ or $D_{n}$ to semi-standard tableaux.

\noindent We define an order on the vertices of the above crystal graphs by
setting%
\begin{gather*}
\mathcal{A}_{n}^{B}=\{\overline{n}<\cdot\cdot\cdot<\overline{1}<0<1<\cdot
\cdot\cdot<n\}\\
\mathcal{A}_{n}^{C}=\{\overline{n}<\cdot\cdot\cdot<\overline{1}<1<\cdot
\cdot\cdot<n\}\text{ and }\\
\mathcal{A}_{n}^{D}=\{\overline{n}<\cdot\cdot\cdot<\overline{2}<%
\begin{tabular}
[c]{l}%
$1$\\
$\overline{1}$%
\end{tabular}
<2<\cdot\cdot\cdot<n\}.
\end{gather*}
Note that $\mathcal{A}_{n}^{D}$ is only partially ordered. For any letter $x$
we set $\overline{\overline{x}}=x.$ Our convention for labelling the crystal
graph of the vector representations are not those used by Kashiwara and
Nakashima.\ To obtain the original description of $B(\lambda)$ from that used
in the sequel it suffices to change each letter $k\in\{1,...,n\}$ into
$\overline{n-k+1}$ and each letter $\overline{k}\in\{\overline{1}%
,...,\overline{n}\}$ into $n-k+1.\;$The interest of this change of convention
is to yield a natural extension of the above alphabets$.$

\bigskip

For types $B_{n},C_{n}$ and $D_{n},$ we identify the vertices of the crystal
graph $G_{n}^{B}=\underset{l}{%
{\textstyle\bigoplus}
}B(\Lambda_{n-1}^{B})^{\bigotimes l},G_{n}^{C}=\underset{l}{%
{\textstyle\bigoplus}
}B(\Lambda_{n-1}^{C})^{\bigotimes l}$ and $G_{n}^{D}=\underset{l}{%
{\textstyle\bigoplus}
}B(\Lambda_{n-1}^{D})^{\bigotimes l}$ respectively with the words on
$\mathcal{A}_{n}^{B},\mathcal{A}_{n}^{C}$ and $\mathcal{A}_{n}^{D}$.$\;$For
any $w\in G_{n}$ we have $\mathrm{wt}(w)=d_{\overline{n}}\varepsilon
_{\overline{n}}+d_{\overline{n-1}}\varepsilon_{\overline{n-1}}\cdot\cdot
\cdot+d_{\overline{1}}\varepsilon_{\overline{1}}$ where for all $i=1,...,n$
$d_{\overline{i}}$ is the number of letters $\overline{i}$ of $w$ minus its
number of letters $i.$

\noindent Consider $\lambda$ a generalized partition with nonnegative parts.
Suppose first that $\lambda$ is a partition.\ Write $T_{\lambda}$ for the
filling of the Young diagram of shape $\lambda$ whose $k$-th row contains only
letters $\overline{n-k+1}.$ Let $b_{\lambda}$ be the vertex of $B(\Lambda
_{n-1})^{\bigotimes\left|  \lambda\right|  }$ obtained by column reading
$T_{\lambda}$ from right to left and top to bottom$.$ Kashiwara and Nakashima
realize $B(\lambda)$ as the connected component of the tensor power
$B(\Lambda_{n-1})^{\bigotimes\left|  \lambda\right|  }$ of highest weight
vertex $b_{\lambda}$. For each roots system $B_{n},C_{n}$ and $D_{n},$ the
Kashiwara-Nakashima tableaux of type $B_{n},C_{n}$, $D_{n}$ and shape
$\lambda$ are defined as the tableaux whose column readings are the vertices
of $B(\lambda)$. We will denote by $\mathrm{w}(T)$ the column reading of the
tableau $T.$

\noindent Now suppose that $\lambda$ belongs to $P_{+}^{B_{n}}$ or
$P_{+}^{D_{n}}$ and its parts are half nonnegative integers. In this case we
can write $\lambda=\lambda%
{{}^\circ}%
+(1/2,...,1/2)$ with $\lambda%
{{}^\circ}%
$ a partition and $B(\lambda)$ can be realized as the connected component of
the crystal graph $\frak{G}_{n}^{0}=B(\Lambda_{n-1})^{\bigotimes\left|
\lambda%
{{}^\circ}%
\right|  }\otimes B(\Lambda_{0})$ of highest weight vertex $b_{\lambda
}=b_{\lambda%
{{}^\circ}%
}\otimes b_{\Lambda_{0}}$ where $b_{\Lambda_{0}}$ is the highest weight vertex
of $B(\Lambda_{0})$ the crystal graph of the spin representation
$V(\Lambda_{0})$ of the corresponding Lie algebra$.$ The vertices of
$B(\Lambda_{0})$ are labelled by spin columns which are special column shaped
diagrams of width $1/2$ and height $n.$ Then the vertices of $B(\lambda)$ can
be identified with the column readings of Kashiwara-Nakashima's spin tableaux
of types $B_{n},D_{n}$ and shape $\lambda$ obtained by adding a column shape
diagram of width $1/2$ to the Young diagram associated to $\lambda%
{{}^\circ}%
.$

\noindent Finally suppose that $\lambda$ belongs to $P_{+}^{D_{n}}$ and
verifies $\lambda_{\overline{1}}<0.$ The above description of $B(\lambda)$
remain valuable up to the following minor modifications. If the parts of
$\lambda$ are integers the letters $\overline{1}$ must be changed into letters
$1$ in the above definition of $T_{\lambda}.$ Otherwise we set $\lambda
=\lambda%
{{}^\circ}%
+(1/2,...,1/2,-1/2)$ where $\lambda%
{{}^\circ}%
$ is generalized partition with integer parts.\ Then $B(\lambda)$ is realized
as the connected component of the crystal graph $\frak{G}_{n}^{1}%
=B(\Lambda_{n-1})^{\bigotimes\left|  \lambda%
{{}^\circ}%
\right|  }\otimes B(\Lambda_{1}^{D_{n}})$ of highest weight vertex
$b_{\lambda}=b_{\lambda^{%
{{}^\circ}%
}}\otimes b_{\Lambda_{1}^{D_{n}}}$ where $b_{\Lambda_{1}^{D_{n}}}$ is the
highest weight vertex of $B(\Lambda_{1}^{D_{n}})$ the crystal graph of the
spin representation $V(\Lambda_{1}^{D_{n}})$.

\bigskip

For any generalized partition $\lambda$ of length $n,$ write $\mathbf{T}%
^{B_{n}}(\lambda),$ $\mathbf{T}^{C_{n}}(\lambda)$ and $\mathbf{T}^{D_{n}%
}(\lambda)$ respectively for the sets of Kashiwara-Nakashima's tableaux of
shape $\lambda$. Set $\mathbf{T}^{B_{n}}=\underset{\lambda\in P_{B_{n}}^{+}%
}{\cup}\mathbf{T}^{B_{n}}(\lambda),$ $\mathbf{T}^{C_{n}}=\underset{\lambda\in
P_{C_{n}}^{+}}{\cup}\mathbf{T}^{C_{n}}(\lambda)$ and $\mathbf{T}^{D_{n}%
}=\underset{\lambda\in P_{D_{n}}^{+}}{\cup}\mathbf{T}^{D_{n}}(\lambda).$ In
the sequel we only summarize the combinatorial description of the partition
shaped tableaux that is, tableaux of $\mathbf{T}^{n}(\lambda)$ where the parts
of $\lambda$ are integers (with eventually $\lambda_{\overline{1}}<0$ for the
root system $D_{n}).$ We refer the reader to \cite{Ba}, \cite{KN}, \cite{Lec}
and \cite{lec2} for the complete description of $\mathbf{T}^{n}(\lambda)$
which necessitates a large amount of combinatorial definitions especially when
the parts of $\lambda$ are half integers.

\noindent So consider $\lambda$ a generalized partition with integer parts.
Suppose first that $\lambda_{\overline{n}}=1.$ Then the tableaux of
$\mathbf{T}^{n}(\lambda)$ are called the $n$-admissible columns. The
$n$-admissible columns of types $B_{n},C_{n}$ and $D_{n}$ are in particular
columns of types $B_{n},C_{n}$ and $D_{n}$ that is have the form%
\begin{equation}
C=%
\begin{tabular}
[c]{|l|}\hline
$C_{-}$\\\hline
$C_{0}$\\\hline
$C_{+}$\\\hline
\end{tabular}
,C=%
\begin{tabular}
[c]{|l|}\hline
$C_{-}$\\\hline
$C_{+}$\\\hline
\end{tabular}
\text{ and }C=%
\begin{tabular}
[c]{|c|}\hline
$D_{-}$\\\hline
$D$\\\hline
$D_{+}$\\\hline
\end{tabular}
\label{col}%
\end{equation}
where $C_{-},C_{+},C_{0},D_{-},D_{+}$ and $D$ are column shaped Young diagrams
such that%
\[
\left\{
\begin{tabular}
[c]{l}%
$C_{-}$ is filled by strictly increasing barred letters from top to bottom\\
$C_{+}$ is filled by strictly increasing unbarred letters from top to bottom\\
$C_{0}$ is filled by letters $0$\\
$D_{-}$ is filled by strictly increasing letters $\leq\overline{2}$ from top
to bottom\\
$D_{+}$ is filled by strictly increasing letters $\geq2$ from top to bottom\\
$D$ is filled by letters $\overline{1}$ or $1$ with differents letters in two
adjacent boxes
\end{tabular}
\right.  .
\]
Note that all the columns are not $n$-admissible even if their letters $a$
satisfy $\overline{n}\leq a\leq n.$ More precisely a column $C$ of (\ref{col})
is $n$-admissible if and only if it can be duplicated following a simple
algorithm described in \cite{lec2} into a pair $(lC,rC)$ of columns without
pair of opposite letters $(x,\overline{x})$ (the letter $0$ is counted as the
pair $(0,\overline{0})$) and containing only letters $a$ such that
$\overline{n}\leq a\leq n.$

\begin{example}
For the column $C=%
\begin{tabular}
[c]{|l|}\hline
$\mathtt{\bar{3}}$\\\hline
$\mathtt{\bar{1}}$\\\hline
$\mathtt{0}$\\\hline
$\mathtt{1}$\\\hline
$\mathtt{2}$\\\hline
\end{tabular}
$ of type $B$ we have $lC=$%
\begin{tabular}
[c]{|l|}\hline
$\mathtt{\bar{5}}$\\\hline
$\mathtt{\bar{4}}$\\\hline
$\mathtt{\bar{3}}$\\\hline
$\mathtt{1}$\\\hline
$\mathtt{2}$\\\hline
\end{tabular}
and $rC=$%
\begin{tabular}
[c]{|l|}\hline
$\mathtt{\bar{3}}$\\\hline
$\mathtt{\bar{1}}$\\\hline
$\mathtt{3}$\\\hline
$\mathtt{4}$\\\hline
$\mathtt{5}$\\\hline
\end{tabular}
. Hence $C$ is $5$-admissible but not $n$-admissible for $n\leq4.$
\end{example}

\noindent Now for a general $\lambda$ with integer parts, a tableau
$T\in\mathbf{T}^{n}(\lambda)$ can be regarded as a filling of the Young
diagram of shape $\lambda$ if $\lambda_{\overline{1}}\geq0$ (of shape
$\lambda^{\ast}$ otherwise) such that

\begin{itemize}
\item $T=C_{1}\cdot\cdot\cdot C_{r}$ where the columns $C_{i}$ of $T$ are $n$-admissible,

\item  for any $i\in\{1,...r-1\}$ the columns of the tableau $r(C_{i}%
)l(C_{i+1})$ weakly increase from left to right and do not contain special
configurations (detailed in \cite{KN} and \cite{lec2}) when $T$ is of type
$D_{n}.$
\end{itemize}

\noindent\textbf{Remark: }We have $\mathbf{T}^{n}(\lambda)\subset
\mathbf{T}^{n+1}(\lambda^{\#})$ where $\lambda^{\#}=(\lambda_{\overline{n}%
},...,\lambda_{\overline{1}},0)$ since the $n$-admissible columns are also
$(n+1)$-admissible and the duplication process of a column does not depend on
$n.$ To simplify the notation we will write in the sequel $\mathbf{T}%
^{n+1}(\lambda)$ instead of $\mathbf{T}^{n+1}(\lambda^{\#})$ for any
$\lambda\in P_{n}^{+}.$

\subsection{Insertion schemes and plactic monoids}

There exist insertion schemes related to each classical root system \cite{Ba},
\cite{Lec} and \cite{lec2} analogous for Kashiwara-Nakashima's tableaux to the
well known bumping algorithm on semi-standard tableaux.

\noindent Denote by $\sim_{n}^{B},\sim_{n}^{C}$ and $\sim_{n}^{D}$ the
equivalence relations defining on the vertices of $G_{n}^{B},G_{n}^{C}$ and
$G_{n}^{D}$ by $w_{1}\sim_{n}w_{2}$ if and only if $w_{1}$ and $w_{2}$ belong
to the same connected component of $G_{n}.\;$For any word $w,$ the insertions
schemes permit to compute the unique tableau $P_{n}(w)$ such that $w\sim
_{n}\mathrm{w}(P_{n}(w))$. In fact $\sim_{n}^{B},\sim_{n}^{C}$ and $\sim
_{n}^{D}$ are congruencies $\equiv_{n}^{B},\equiv_{n}^{C}$ and $\equiv_{n}%
^{D}$ \cite{Ba} \cite{Lec} \cite{lec2} \cite{lit} obtained respectively as the
quotient of the free monoids of words on $\mathcal{A}_{n}^{C},\mathcal{A}%
_{n}^{B}$ and $\mathcal{A}_{n}^{D}$ by two kinds of relations.

\noindent The first is constituted by relations of length $3$ analogous to
Knuth relations defining Lascoux-Sch\"{u}tzenberger's plactic monoid. In fact
these relations are precisely those which are needed to describe the insertion
$x\rightarrow C$ of a letter $x$ in a $n$-admissible column $C=%
\begin{tabular}
[c]{|l|}\hline
$a$\\\hline
$b$\\\hline
\end{tabular}
$ such that
\begin{tabular}
[c]{|l|}\hline
$a$\\\hline
$b$\\\hline
$x$\\\hline
\end{tabular}
is not a column. This can be written%
\begin{equation}
x\rightarrow%
\begin{tabular}
[c]{|l|}\hline
$a$\\\hline
$b$\\\hline
\end{tabular}
=%
\begin{tabular}
[c]{c|c|}\cline{2-2}%
& $a$\\\hline
\multicolumn{1}{|c|}{$x$} & $b$\\\hline
\end{tabular}
=%
\begin{tabular}
[c]{|l|l}\hline
$a^{\prime}$ & \multicolumn{1}{|l|}{$x^{\prime}$}\\\hline
$b^{\prime}$ & \\\cline{1-1}%
\end{tabular}
\label{trans_ele}%
\end{equation}
and contrary to the insertion scheme for the semi-standard tableaux the sets
$\{a^{\prime},b^{\prime},x^{\prime}\}$ and $\{a,b,c\}$ are not necessarily
equal (i.e. the relations are not homogeneous in general).

\noindent Next we have the contraction relations which do not preserve the
length of the words.\ These relations are precisely those which are needed to
describe the insertion $x\rightarrow C$ of a letter $x$ such that
$\overline{n}\leq x\leq n$ in a $n$-admissible column $C$ such that
\begin{tabular}
[c]{|l|}\hline
$C$\\\hline
$x$\\\hline
\end{tabular}
(obtained by adding the letter $x$ on bottom of $C)$ is a column which is not
$n$-admissible. In this case
\begin{tabular}
[c]{|l|}\hline
$C$\\\hline
$x$\\\hline
\end{tabular}
is necessarily $(n+1)$-admissible and have to be contracted to give a
$n$-admissible column. We obtain $x\rightarrow C=\widetilde{C}$ with
$\widetilde{C}$ a $n$-admissible column of height $h(C)$ or $h(C)-1.$

\noindent The insertion of the letter $x$ in a $n$-admissible column $C$ of
arbitrary height such that
\begin{tabular}
[c]{|l|}\hline
$C$\\\hline
$x$\\\hline
\end{tabular}
is not a column can then be pictured by%
\[
x\rightarrow%
\begin{tabular}
[c]{|c|}\hline
$a_{1}$\\\hline
$\cdot$\\\hline
$a_{k-2}$\\\hline
$a_{k-1}$\\\hline
$a_{k}$\\\hline
\end{tabular}
=%
\begin{tabular}
[c]{c|c|}\cline{2-2}%
\ \ \ \ \  & $a_{1}$\\\cline{2-2}%
& $\cdot$\\\cline{2-2}%
& $a_{k-2}$\\\cline{2-2}%
& $a_{k-1}$\\\hline
\multicolumn{1}{|c|}{$x$} & $a_{k}$\\\hline
\end{tabular}
=%
\begin{tabular}
[c]{c|c}\cline{2-2}%
& \multicolumn{1}{|c|}{$a_{1}$}\\\cline{2-2}%
& \multicolumn{1}{|c|}{$\cdot$}\\\cline{2-2}%
& \multicolumn{1}{|c|}{$a_{k-2}$}\\\hline
\multicolumn{1}{|c|}{$\delta_{k-1}$} & \multicolumn{1}{|c|}{$y$}\\\hline
\multicolumn{1}{|c|}{$d_{k}$} & \\\cline{1-1}%
\end{tabular}
=\cdot\cdot\cdot=%
\begin{tabular}
[c]{|c|c}\hline
$d_{1}$ & \multicolumn{1}{|c|}{$\ \ z$ \ \ }\\\hline
$\cdot$ & \\\cline{1-1}%
$\cdot$ & \\\cline{1-1}%
$d_{k-1}$ & \\\cline{1-1}%
$d_{k}$ & \\\cline{1-1}%
\end{tabular}
\]
that is, one elementary transformation (\ref{trans_ele}) is applied to each
step. One proves that $x\rightarrow C$ is then a tableau of $\mathbf{T}^{(n)}$
with two columns respectively of height $h(C)$ and $1$.

\bigskip

\noindent Now we can define the insertion $x\rightarrow T$ of the letter $x$
such that $\overline{n}\leq x\leq n$ in the tableau $T\in\mathbf{T}%
^{n}(\lambda)$.\ Set $T=C_{1}\cdot\cdot\cdot C_{r}$ where $C_{i},$ $i=1,...,r$
are the $n$-admissible columns of $T.\;$

\begin{enumerate}
\item  When
\begin{tabular}
[c]{|l|}\hline
$C_{1}$\\\hline
$x$\\\hline
\end{tabular}
is not a column, write $x\rightarrow C=%
\begin{tabular}
[c]{|l|l|}\hline
$C_{1}^{\prime}$ & $y$\\\hline
\end{tabular}
$ where $C_{1}^{\prime}$ is an admissible column of height $h(C_{1})$ and $y$
a letter.\ Then $x\rightarrow T=C_{1}^{\prime}(y\rightarrow C_{2}\cdot
\cdot\cdot C_{r})$ that is, $x\rightarrow T$ is the juxtaposition of
$C_{1}^{\prime}$ with the tableau $\widehat{T}$ obtained by inserting $y$ in
the tableau $C_{2}\cdot\cdot\cdot C_{r}.$

\item  When
\begin{tabular}
[c]{|l|}\hline
$C_{1}$\\\hline
$x$\\\hline
\end{tabular}
is a $n$-admissible column, $x\rightarrow T$ is the tableau obtained by adding
a box containing $x$ on bottom of $C_{1}$.\ 

\item  When
\begin{tabular}
[c]{|l|}\hline
$C_{1}$\\\hline
$x$\\\hline
\end{tabular}
is a column which is not $n$-admissible, write $x\rightarrow C=\widetilde{C}$
and set $\mathrm{w}(\widetilde{C})=y_{1}\cdot\cdot\cdot y_{s}$ where the
$y_{i}$'s are letters. Then $x\rightarrow T=y_{s}\rightarrow(y_{s-1}%
\rightarrow(\cdot\cdot\cdot y_{1}\rightarrow\widehat{T}))$ that is
$x\rightarrow T$ is obtained by inserting successively the letters of
$\widetilde{C}$ into the tableau $\widehat{T}=C_{2}\cdot\cdot\cdot C_{r}$.
Note that there is no new contraction during this $s$ insertions.
\end{enumerate}

\noindent\textbf{Remarks:}

\noindent$\mathrm{(i)}\mathbf{:}$ The $P_{n}$-symbol defined above can be
computed recursively by setting $P_{n}(w)=%
\begin{tabular}
[c]{|l|}\hline
$w$\\\hline
\end{tabular}
$ if $w$ is a letter and $P_{n}(w)=x\rightarrow P_{n}(u)$ where $w=ux$ with
$u$ a word and $x$ a letter otherwise.

\noindent$\mathrm{(ii)}\mathbf{:}$ Consider $T\in\mathbf{T}^{n}(\lambda
)\subset\mathbf{T}^{n+1}(\lambda)$ and a letter $x$ such that $\overline
{n}\leq x\leq n.$ The tableau obtained by inserting $x$ in $T$ may depend
wether $T$ is regarded as a tableau of $\mathbf{T}^{n}(\lambda)$ or as a
tableau of $\mathbf{T}^{n+1}(\lambda)$.\ Indeed if \
\begin{tabular}
[c]{|l|}\hline
$C_{1}$\\\hline
$x$\\\hline
\end{tabular}
is not $n$-admissible then it is necessarily $(n+1)$-admissible since $C_{1}$
is $n$-admissible. Hence there is no contraction during the insertion
$x\rightarrow T$ when it is regarded as a tableau of $\mathbf{T}^{n+1}(\lambda).$

\noindent$\mathrm{(iii)}\mathbf{:}$ Consider $w\in\mathcal{A}_{n}$, from
$\mathrm{(ii)}\mathbf{\ }$we deduce that there exists an integer $m\geq n$
minimal such that $P_{m}(w)$ can be computed without using contraction
relation. Then for any $k\geq m,$ $P_{k}(w)=P_{m}(w).$

\noindent$\mathrm{(iv)}\mathbf{:}$ Similarly to the bumping algorithm for
semi-standard tableaux, the insertion algorithms described above are reversible.

\noindent$\mathrm{(v)}\mathbf{:}$ There also exit insertion algorithms for the
spin tableaux of types $B_{n}$ and $D_{n}$ \cite{lec2}. To make the paper more
readable we only establish the combinatorial results contained in the sequel
for the partition shaped tableaux.\ Nevertheless note that they can be
extended to take also into account the spin tableaux associated to the root
systems $B_{n}$ and $D_{n}.$

\begin{lemma}
\label{lem_fact_row}Consider $\lambda,\mu\in P_{n}^{+}.$ Let $T\in
\mathbf{T}^{n}(\lambda).$ If $\lambda$ and $\mu$ have integer parts, then
there exists a unique pair $(R,T^{\prime})$ such that%
\[
\mathrm{w}(T)\equiv_{n}\mathrm{w}(R)\otimes\mathrm{w}(T^{\prime})
\]
where $R\in\mathbf{T}^{n}(\lambda_{\overline{n}}\Lambda_{n-1})$ is a row
tableau of length $\lambda_{\overline{n}}$ and $T^{\prime}\in\mathbf{T}%
^{n-1}(\lambda^{\prime})$ with $\lambda^{\prime}=(\lambda_{\overline{n-1}%
},...,\lambda_{\overline{1}}).$
\end{lemma}

\begin{proof}
When $\lambda_{\overline{1}}\geq0$ we have
\[
b_{\lambda}\equiv_{n}(\overline{n})^{\otimes\lambda_{\overline{n}}}%
\otimes\left(  (\overline{1})^{\otimes\widehat{\lambda}_{\overline{1}}%
^{\prime}}\otimes(\overline{1}\ \overline{2})^{\otimes\widehat{\lambda
}_{\overline{2}}^{\prime}}\otimes\cdot\cdot\cdot\otimes(\overline
{1}\ \overline{2}\cdot\cdot\cdot\overline{n-1})^{\otimes\widehat{\lambda
}_{\overline{n-1}}^{\prime}}\right)  \equiv_{n}b_{\lambda_{\overline{n}%
}\Lambda_{n-1}}\otimes b_{\lambda^{\prime}}%
\]
with the notation used in \ref{subsec_KF}. Indeed the plactic relations on
words containing only barred letters coincide with Knuth relations. This
implies the existence of the required pair $(R,T^{\prime})$. Now if
$\mathrm{w}(T)\equiv_{n}\mathrm{w}(R)\otimes\mathrm{w}(T^{\prime}),$ we deduce
from \ref{lem_plu_hp} that the highest weight vertex of the connected
component of $G_{n}$ containing $\mathrm{w}(R)\otimes\mathrm{w}(T^{\prime})$
is necessarily $b_{\lambda_{\overline{n}}\Lambda_{n-1}}\otimes b_{\lambda
^{\prime}}$. Thus the pair $(R,T^{\prime})$ is unique.
\end{proof}

\noindent\textbf{Remark: }The pair $(R,T^{\prime})$ can be explicitly computed
by using the reverse insertion schemes.

\section{\label{sec_morris}Morris type recurrence formulas for the orthogonal
root systems}

In this section we introduce recurrence formulas for computing Kostka-Foulkes
polynomials analogous for types $B_{n}$ and $D_{n}$ to Morris recurrence
formula.\ They allow to explain the Kostka-Foulkes polynomials for types
$B_{n}$ and $D_{n}$ respectively as combinations of Kostka-Foulkes polynomials
for types $B_{n-1}$ and $D_{n-1}.$ We essentially proceed as we have done in
\cite{lec3} for the root system $C_{n}$. So we only sketch the arguments
except for Theorems \ref{Th_rec_morrB} and \ref{Th_rec_morrD} for which the
proofs necessitate refinements of the proof of Theorem 3.2.1 of \cite{lec3}.

\noindent We classically realize $so_{2n-1},sp_{2n-2}$ and $so_{2n-2}$
respectively as the sub-algebras of $so_{2n+1},sp_{2n}$ and $so_{2n}$
generated by the Chevalley operators $e_{i},f_{i}$ and $t_{i},$ $i=0,...n-2.$
The weight lattice $P_{n-1}$ of these algebras of rank $n-1$ is the
$\mathbb{Z}$-lattice generated by the $\varepsilon_{\overline{i}},$
$i=1,...,n-1$ and $P_{n-1}^{+}=P_{n}^{+}\cap P_{n-1}$ is the set of dominant
weights. The Weyl group $W_{n-1}$ is the sub-group of $W_{n}$ generated by the
$s_{i},$ $i=0,...n-2$ and we have $R_{n-1}^{+}=R_{n}\cap P_{n-1}^{+}.$

\noindent Given any positive integer $r,$ set $B^{B_{n}}(r)=B(r\Lambda
_{n-1}^{B_{n}}),B^{C_{n}}(r)=B(r\Lambda_{n-1}^{C_{n}}),$ and $B^{D_{n}%
}(r)=B(r\Lambda_{n-1}^{D_{n}})$. To obtain our recurrence formulas we need to
describe the decomposition $B(\gamma)\otimes B(r)$ with $\gamma\in P_{n}^{+}$
and $r\geq0$ an integer into its irreducible components.\ This is analogous
for types $B_{n}$ and $D_{n}$ to Pieri rule.

\subsection{Pieri rule for types $B_{n}$ and $D_{n}$}

It follows from \cite{KN} that the vertices of $B^{B_{n}}(r),$ $B^{C_{n}}(r)$
and $B^{D_{n}}(r)$ can be respectively identified to the words
\begin{equation}
L=(n)^{k_{n}}\cdot\cdot\cdot(2)^{k_{2}}(1)^{k_{1}}(\overline{1})^{k_{\bar{1}}%
}(\overline{2})^{k_{\bar{2}}}\cdot\cdot\cdot(\overline{n})^{k_{\bar{n}}%
}\text{, }L=(n)^{k_{n}}\cdot\cdot\cdot(2)^{k_{2}}(1)^{k_{1}}(0)(\overline
{1})^{k_{\bar{1}}}(\overline{2})^{k_{\bar{2}}}\cdot\cdot\cdot(\overline
{n})^{k_{\bar{n}}} \label{LB}%
\end{equation}%
\begin{equation}
L=(n)^{k_{n}}\cdot\cdot\cdot(2)^{k_{2}}(1)^{k_{1}}(\overline{1})^{k_{\bar{1}}%
}(\overline{2})^{k_{\bar{2}}}\cdot\cdot\cdot(\overline{n})^{k_{\bar{n}}}
\label{LC}%
\end{equation}
and
\begin{equation}
L=(n)^{k_{n}}\cdot\cdot\cdot(2)^{k_{2}}(\overline{1})^{k_{\overline{1}}%
}(\overline{2})^{k_{\bar{2}}}\cdot\cdot\cdot(\overline{n})^{k_{\bar{n}}%
}\text{, }L=(n)^{k_{n}}\cdot\cdot\cdot(2)^{k_{2}}(1)^{k_{1}}(\overline
{2})^{k_{\bar{2}}}\cdot\cdot\cdot(\overline{n})^{k_{\bar{n}}} \label{LD}%
\end{equation}
of length $r$ where $k_{\overline{i}},k_{i}$ are positive integers, $(x)^{k}$
means that the letter $x$ is repeated $k$ times in $L.$ Note that there can be
only one letter $0$ in the vertices of $B^{B_{n}}(r)$ and the letters
$\overline{1}$ and $1$ can not appear simultaneously in the vertices of
$B^{D_{n}}(r).$

Let $\gamma=(\gamma_{\overline{n}},...,\gamma_{\overline{1}})\in P_{n}^{+}.$
When $\gamma\in P_{B_{n}}^{+}$ set $B(\gamma)\otimes B^{B_{n}}(r)=\underset
{\lambda\in P_{B_{n}}^{+}}{\cup}B(\lambda)^{b_{\gamma,r}^{\lambda}}$ that is
$b_{\gamma,r}^{\lambda}$ is the multiplicity of $V(\lambda)$ in $V(\gamma
)\otimes V(r\Lambda_{n-1}^{B}).$ Similarly set $B(\gamma)\otimes B^{C_{n}%
}(r)=\underset{\lambda\in P_{C_{n}}^{+}}{\cup}B(\lambda)^{c_{\gamma
,r}^{\lambda}}$ and $B(\gamma)\otimes B^{D_{n}}(r)=\underset{\lambda\in
P_{D_{n}}^{+}}{\cup}B(\lambda)^{d_{\gamma,r}^{\lambda}}$ when $\gamma$ belongs
respectively to $P_{C_{n}}^{+}$ and $P_{D_{n}}^{+}.$

\noindent Write $b_{\gamma}$ for the highest weight vertex of $B(\gamma).$ The
two following lemmas and their corollaries are consequences of Lemma
\ref{lem_plu_hp}.

\begin{lemma}
\label{lem_b_gamm_tens_Lb}$b_{\gamma}\otimes L$ is a highest weight vertex of
$B(\gamma)\otimes B^{B_{n}}(r)$ if and only if the following conditions holds:

\noindent$\mathrm{(i):}$ $\gamma_{\overline{1}}-k_{1}\geq0$ if $k_{0}=0$,
$\gamma\overline{_{1}}-k_{1}>0$ otherwise

\noindent$\mathrm{(ii):}$ $\gamma_{\overline{i+1}}-k_{i+1}\geq\gamma
_{\overline{i}}$ for $i=1,...,n-1$

\noindent$\mathrm{(iii):}$ $\gamma_{\overline{i}}-k_{i}+k_{\overline{i}}%
\leq\gamma_{\overline{i+1}}-k_{i+1}$ for $i=1,...,n-1$
\end{lemma}

\begin{corollary}
\label{cor_pieriB}The multiplicity $b_{\gamma,r}^{\lambda}$ is the number of
vertices $L\in B^{B_{n}}(r)$ such that $k_{\overline{i}}-k_{i}=\lambda
_{\overline{i}}-\gamma_{\overline{i}}$ for $i=1,....,n$ and

\noindent$\mathrm{(i)}:$ $\lambda_{\overline{i}}\leq\lambda_{\overline{i+1}%
}-k_{\overline{i+1}}$ for $i=1,...,n-1$,

\noindent$\mathrm{(ii)}:$ $\lambda_{\overline{i+1}}-k_{\overline{i+1}}%
\geq\lambda_{\overline{i}}+k_{i}-k_{\overline{i}}$ for $i=1,...,n-1$,

\noindent$\mathrm{(iii)}:$ $\lambda_{\overline{1}}-k_{\overline{1}}\geq0$ if
$k_{0}=0$ (i.e. $k_{\overline{1}}+\cdot\cdot\cdot+k_{\overline{n}}+k_{1}%
+\cdot\cdot\cdot+k_{n}=r)$ and $\lambda_{\overline{1}}-k_{\overline{1}}>0$
otherwise (i.e. $k_{\overline{1}}+\cdot\cdot\cdot+k_{\overline{n}}+k_{1}%
+\cdot\cdot\cdot+k_{n}=r-1$).
\end{corollary}

\begin{lemma}
\label{lem_b_gamm_tens_Ld}$b_{\gamma}\otimes L$ is a highest weight vertex of
$B(\gamma)\otimes B^{D_{n}}(r)$ if and only if the following conditions holds:

\noindent$\mathrm{(i):}$ $\gamma_{\overline{2}}-k_{2}\geq\gamma_{\overline{1}%
}$ if $\gamma_{\overline{1}}\geq0$ and $\gamma_{\overline{2}}-k_{2}\geq
-\gamma_{\overline{1}}$ otherwise,

\noindent$\mathrm{(ii):}$ $\gamma_{\overline{i+1}}-k_{i+1}\geq\gamma
_{\overline{i}}$ for $i=2,...,n-1$,

\noindent$\mathrm{(iii):}\gamma_{\overline{1}}+k_{\overline{1}}\leq
\gamma_{\overline{2}}-k_{2}$ if $k_{1}=0$ and $-\gamma_{\overline{1}}%
+k_{1}\leq\gamma_{\overline{2}}-k_{2}$ otherwise,

\noindent$\mathrm{(iii):}$ $\gamma_{\overline{i}}-k_{i}+k_{\overline{i}}%
\leq\gamma_{\overline{i+1}}-k_{i+1}$ for $i=2,...,n-1$
\end{lemma}

\begin{corollary}
\label{cor_pieriD}The multiplicity $d_{\gamma,r}^{\lambda}$ is the number of
vertices $L\in\otimes B^{D_{n}}(r)$ such that $k_{\overline{i}}-k_{i}%
=\lambda_{\overline{i}}-\gamma_{\overline{i}}$ for $i=1,....,n$, and

\noindent$\mathrm{(i)}:$ $\lambda_{\overline{1}}\leq\lambda_{\overline{2}%
}-k_{\overline{2}}$ if $k_{1}=0$ and $-\lambda_{\overline{1}}\leq
\lambda_{\overline{2}}-k_{\overline{2}}$ otherwise

\noindent$\mathrm{(ii)}:\lambda_{\overline{i}}\leq\lambda_{\overline{i+1}%
}-k_{\overline{i+1}}$ for $i=2,...,n-1$,

\noindent$\mathrm{(iii)}:$ $\lambda_{\overline{i+1}}-k_{\overline{i+1}}%
\geq\lambda_{\overline{i}}+k_{i}-k_{\overline{i}}$ for $i=2,...,n-1$,

\noindent$\mathrm{(iv)}:$ $\left\{
\begin{tabular}
[c]{l}%
$(a):\lambda_{\overline{2}}-k_{\overline{2}}\geq\lambda_{\overline{1}%
}-k_{\overline{1}}$ if $k_{1}=0$ and $\gamma_{\overline{1}}\geq0$\\
$(b):\lambda_{\overline{2}}-k_{\overline{2}}\geq\lambda_{\overline{1}}+k_{1}$
if $k_{\overline{1}}=0$ and $\gamma_{\overline{1}}\geq0$%
\end{tabular}
\right\}  $, $\left\{
\begin{tabular}
[c]{l}%
$(c):\lambda_{\overline{2}}-k_{\overline{2}}\geq-\lambda_{\overline{1}}-k_{1}$
if $k_{\overline{1}}=0$ and $\gamma_{\overline{1}}<0$\\
$(d):\lambda_{\overline{2}}-k_{\overline{2}}\geq-\lambda_{\overline{1}%
}+k_{\overline{1}}$ if $k_{1}=0$ and $\gamma_{\overline{1}}<0$%
\end{tabular}
\right\}  .$
\end{corollary}

\noindent\textbf{Remarks: }

\noindent$\mathrm{(i):}$\textbf{ }In the above corollaries, $b_{\gamma
,r}^{\lambda}$ and $d_{\gamma,r}^{\lambda}$ are the number of ways of starting
with $\gamma,$ removing a horizontal strip to obtain a partition $\nu$
(corresponding to the unbarred letters of $L)$ and then adding a horizontal
strip (corresponding to the barred letters of $L$) to obtain $\lambda.$

\noindent$\mathrm{(ii):}$ $B(\gamma)\otimes B((r)_{n})$ is not multiplicity
free in general.

\noindent$\mathrm{(iii):}$ Consider $\gamma=(\gamma_{\overline{n}}%
,...,\gamma_{\overline{1}})\in P_{B_{n}}$ (resp. $P_{D_{n}})$ such that
$\lambda=(\lambda_{\overline{n}},...,\lambda_{\overline{1}})\in P_{B_{n}}$
(resp. $P_{D_{n}})$ defined by $\lambda_{\overline{i}}=\gamma_{\overline{i}%
}+k_{\overline{i}}-k_{\overline{i}},$ $i=1,...,n$ verifies conditions
$\mathrm{(i),(ii)}$ and $\mathrm{(iii)}$ of Corollary \ref{cor_pieriB} (resp.
\ref{cor_pieriD}). Then $\gamma\in P_{B_{n}}^{+}$ (resp. $P_{D_{n}}^{+})$ that
is $\gamma$ is a generalized partition.

\subsection{Recurrence formulas}

Consider $\gamma\in P_{n}^{+}$ and $r$ a positive integer. We set%
\begin{gather*}
\left(  \gamma\otimes r\right)  _{B_{n}}=\{\lambda\in P_{B_{n}}^{+},\text{
}b_{\gamma,r}^{\lambda}\neq0\},\left(  \gamma\otimes r\right)  _{C_{n}%
}=\{\lambda\in P_{C_{n}}^{+},\text{ }c_{\gamma,r}^{\lambda}\neq0\}\\
\text{and }\left(  \gamma\otimes r\right)  _{D_{n}}=\{\lambda\in P_{D_{n}}%
^{+},\text{ }d_{\gamma,r}^{\lambda}\neq0\}.
\end{gather*}

\noindent For the root system $C_{n}$ and $\mu=(\mu_{\overline{n}}%
,...,\mu_{\overline{1}})$, we have established in \cite{lec3} the following
analogue of Morris recurrence formula:%
\[
Q_{\mu}^{C_{n}}=\sum_{\gamma\in P_{C_{n-1}}^{+}}\sum_{R=0}^{+\infty}%
\sum_{r+2m=R}q^{m+r}\sum_{\lambda\in\left(  \gamma\otimes r\right)  _{C_{n-1}%
}}c_{\gamma,r}^{\lambda}K_{\lambda,\mu^{\prime}}^{C_{n-1}}(q)s_{(\mu
_{\overline{n}}+R,\gamma)}%
\]
where $\mu^{\prime}=(\mu_{\overline{n-1}},...,\mu_{\overline{1}})\in
P_{C_{n-1}}^{+}.$

\begin{theorem}
\label{Th_rec_morrB}Let $\mu\in P_{B_{n}}^{+}.$ Then%
\begin{equation}
Q_{\mu}^{B_{n}}=\sum_{\gamma\in P_{B_{n-1}}^{+}}\sum_{R=0}^{+\infty}%
\sum_{r+2m=R}q^{R}\sum_{\lambda\in\left(  \gamma\otimes r\right)  _{B_{n-1}}%
}b_{\gamma,r}^{\lambda}K_{\lambda,\mu^{\prime}}^{B_{n-1}}(q)s_{(\mu
_{\overline{n}}+R,\gamma)}. \label{rec_mor_b}%
\end{equation}
\end{theorem}

\begin{proof}
From $Q_{\mu}=\left(  \prod_{\alpha\in R_{B_{n}}^{+}}\dfrac{1}{1-qR_{\alpha}%
}\right)  s_{\mu}$ and Proposition 3.5 of \cite{NR} we can write
\[
Q_{\mu}=\left(  \underset{\alpha\notin R_{B_{n-1}}^{+}}{\prod_{\alpha\in
R_{B_{n}}^{+}}}\dfrac{1}{1-qR_{\alpha}}\right)  \left[  \left(  \underset
{\alpha\in R_{B_{n-1}}^{+}}{\prod}\dfrac{1}{1-qR_{\alpha}}\right)  s_{\mu
}\right]  .
\]
Then by applying Theorem \ref{th_hall_kostka}, we obtain%
\begin{equation}
Q_{\mu}=\left(  \underset{\alpha\notin R_{B_{n-1}}^{+}}{\prod_{\alpha\in
R_{B_{n}}^{+}}}\dfrac{1}{1-qR_{\alpha}}\right)  \left(  \sum_{\lambda\in
P_{B_{n-1}}^{+}}K_{\lambda,\mu^{\prime}}^{B_{n-1}}(q)s_{(\mu_{\overline{n}%
},\lambda)}\right)  . \label{for_pro_q_mu}%
\end{equation}
Set $R_{\overline{i}}=R_{\varepsilon_{\overline{n}}-\varepsilon_{\overline{i}%
}}$ for $i=1,...,n-1$ $R_{n}=R_{\varepsilon\overline{_{n}}}$ and
$R_{i}=R_{\varepsilon_{\overline{n}}+\varepsilon_{\overline{i}}}$ for
$i=1,...,n.$ Recall that for any $\beta\in P_{B_{n-1}},$ $R_{\overline{i}%
}(\beta)=\beta+\varepsilon_{\overline{n}}-\varepsilon_{\overline{i}}$ and
$R_{i}(\beta)=\beta+\varepsilon_{\overline{n}}+\varepsilon_{\overline{i}}.$
Then (\ref{for_pro_q_mu}) implies%
\begin{multline*}
Q_{\mu}=\sum_{\lambda\in P_{B_{n-1}}^{+}}K_{\lambda,\mu^{\prime}}^{B_{n-1}%
}(q)\times\\
\left(  \sum_{r=0}^{+\infty}\sum_{b=0}^{+\infty}\sum_{k_{\overline{1}}%
+\cdot\cdot\cdot+k_{\overline{n-1}}+k_{1}+\cdot\cdot\cdot+k_{n-1}=r}%
q^{r+b}(R_{n})^{b}(R_{1})^{k_{1}}(R_{\overline{1}})^{k_{\overline{1}}}%
\cdot\cdot\cdot(R_{n-1})^{k_{n-1}}(R_{\overline{n-1}})^{k_{\overline{n-1}}%
}s_{(\mu_{\overline{n}},\lambda)}\right)  .
\end{multline*}%
\begin{multline*}
Q_{\mu}=\sum_{r=0}^{+\infty}\sum_{b=0}^{+\infty}q^{r+b}\sum_{\lambda\in
P_{B_{n-1}}^{+}}K_{\lambda,\mu^{\prime}}^{B_{n-1}}(q)\sum_{k_{\overline{1}%
}+\cdot\cdot\cdot+k_{\overline{n-1}}+k_{1}+\cdot\cdot\cdot+k_{n-1}=r}%
s_{(\mu_{\overline{n}}+r+b,\lambda_{\overline{n-1}}+k_{n-1}-k_{\overline{n-1}%
},\cdot\cdot\cdot,\lambda_{\overline{1}}+k_{1}-k_{\overline{1}})}=\\
\sum_{R=0}^{+\infty}\sum_{r=0}^{R}q^{R}\sum_{\lambda\in P_{B_{n-1}}^{+}%
}K_{\lambda,\mu^{\prime}}^{B_{n-1}}(q)\sum_{k_{\overline{1}}+\cdot\cdot
\cdot+k_{\overline{n-1}}+k_{1}+\cdot\cdot\cdot+k_{n-1}=r}s_{(\mu_{\overline
{n}}+R,\lambda_{\overline{n-1}}+k_{n-1}-k_{\overline{n-1}},\cdot\cdot
\cdot,\lambda_{\overline{1}}+k_{1}-k_{\overline{1}})}%
\end{multline*}
by setting $R=r+b.$ Now fix $\lambda,R>0$ and $0<r\leq R$ and write
\begin{align*}
S_{1}  &  =\sum_{k_{\overline{1}}+\cdot\cdot\cdot+k_{\overline{n-1}}%
+k_{1}+\cdot\cdot\cdot+k_{n-1}=r}s_{(\mu_{\overline{n}}+R,\lambda
_{\overline{n-1}}+k_{n-1}-k_{\overline{n-1}},\cdot\cdot\cdot,\lambda
_{\overline{1}}+k_{1}-k_{\overline{1}})}\text{,}\\
S_{2}  &  =\sum_{k_{\overline{1}}+\cdot\cdot\cdot+k_{\overline{n-1}}%
+k_{1}+\cdot\cdot\cdot+k_{n-1}=r-1}s_{(\mu_{\overline{n}}+R,\lambda
_{\overline{n-1}}+k_{n-1}-k_{\overline{n-1}},\cdot\cdot\cdot,\lambda
_{\overline{1}}+k_{1}-k_{\overline{1}})},
\end{align*}
$S_{R,r}=S_{1}+S_{2}$ and $\gamma=(\lambda_{\overline{n-1}}+k_{n-1}%
-k_{\overline{n-1}},...,\lambda_{\overline{1}}+k_{1}-k_{\overline{1}}).$

\noindent$\mathrm{(a):}$\textrm{ }Consider\textrm{ }$\gamma$ appearing in
$S_{1}$ or $S_{2}$ and suppose that there exists $i\in\{1,...,n-2\}$ such that
$\lambda_{\overline{i}}>\lambda_{\overline{i+1}}-k_{\overline{i+1}}.$ Set
$\widetilde{\gamma}=s_{i}\circ\gamma$ that is
\[
\widetilde{\gamma}=s_{i}(\gamma_{\overline{n-1}}+n-3/2,...,\gamma
_{\overline{i+1}}+n-i+1/2,\gamma_{\overline{i}}+n-i,...,\gamma_{\overline{1}%
}+1/2)-(n-3/2,...,1/2).
\]
Then $\gamma_{\overline{s}}=\widetilde{\gamma}_{\overline{s}}$ for $s\neq
i+1,i$, $\widetilde{\gamma}_{\overline{i+1}}=\gamma_{\overline{i}}-1$ and
$\widetilde{\gamma}_{\overline{i}}=,\gamma_{\overline{i+1}}+1$ that is
\[
\left\{
\begin{tabular}
[c]{l}%
$\widetilde{\gamma}_{\overline{i+1}}=\lambda_{\overline{i}}+k_{i}%
-k_{\overline{i}}-1$\\
$\widetilde{\gamma}_{\overline{i}}=\lambda_{\overline{i+1}}+k_{i+1}%
-k_{\overline{i+1}}+1$%
\end{tabular}
\right.  .
\]
Write $\widetilde{k}_{i+1}=k_{i}$, $\widetilde{k}_{i}=k_{i+1}$, $\widetilde
{k}_{\overline{i+1}}=\lambda_{\overline{i+1}}-\lambda_{\overline{i}%
}+k_{\overline{i}}+1$ and $\widetilde{k}_{\overline{i}}=\lambda_{\overline{i}%
}-\lambda_{\overline{i+1}}+k_{\overline{i+1}}-1.$ To make our notation
homogeneous set $\widetilde{k}_{t}=k_{t}$ for any $t\neq i,i+1,\overline
{i},\overline{i+1}.$ Then $\lambda_{\overline{i}}>\lambda_{\overline{i+1}%
}-\widetilde{k}_{\overline{i+1}}.$ We have $\widetilde{k}_{\overline{i+1}}%
\geq0$ and $\widetilde{k}_{\overline{i}}=\lambda_{\overline{i}}-\lambda
_{\overline{i+1}}+k_{\overline{i+1}}-1\geq0$ since $\lambda_{\overline{i}%
}>\lambda_{\overline{i+1}}-k_{\overline{i+1}}$.$\;$Moreover $\widetilde
{k}_{\overline{1}}+\cdot\cdot\cdot+\widetilde{k}_{\overline{n-1}}%
+\widetilde{k}_{1}+\cdot\cdot\cdot+\widetilde{k}_{n-1}=k_{\overline{1}}%
+\cdot\cdot\cdot+k_{\overline{n-1}}+k_{1}+\cdot\cdot\cdot+k_{n-1}$ and for any
$s\in\{1,...,n-2\}$%
\[
\widetilde{\gamma}_{\overline{s}}=\lambda_{\overline{s}}+\widetilde{k}%
_{s}-\widetilde{k}_{\overline{s}}.
\]
$\mathrm{(b):}$ Consider\textrm{ }$\gamma$ appearing in $S_{1}$ and suppose
that $\lambda_{\overline{i}}\leq\lambda_{\overline{i+1}}-k_{\overline{i+1}}$
for all $i=1,...,n-2$ and $\lambda_{\overline{1}}-k_{\overline{1}}<0$. Set
$\widetilde{\gamma}=s_{0}\circ\gamma.$ Then $\gamma_{\overline{s}}%
=\widetilde{\gamma}_{\overline{s}}$ for $s\neq1$ and $\widetilde{\gamma
}_{\overline{1}}=-\lambda_{\overline{1}}-k_{1}+k_{\overline{1}}-1$. Write
$\widetilde{k}_{i}=k_{i},$ $\widetilde{k}_{\overline{i}}=k_{\overline{i}}$ for
all $i=2,...,n-1$ and set $\widetilde{k}_{1}=k_{\overline{1}}-\lambda
_{\overline{1}}-1,$ $\widetilde{k}_{\overline{1}}=k_{1}+\lambda_{\overline{1}%
}.$ We have $\widetilde{k}_{1}\geq0$, $\lambda_{\overline{i}}\leq
\lambda_{\overline{i+1}}-\widetilde{k}_{\overline{i+1}}$ for all $i=1,...,n-2$
and $\lambda_{\overline{1}}-\widetilde{k}_{\overline{1}}\leq0.\;$Moreover
$\widetilde{k}_{\overline{1}}+\cdot\cdot\cdot+\widetilde{k}_{\overline{n-1}%
}+\widetilde{k}_{1}+\cdot\cdot\cdot+\widetilde{k}_{n-1}=r-1$ (thus $\gamma$
appears in $S_{2})$ and $\widetilde{\gamma}_{\overline{1}}=\lambda
_{\overline{1}}+\widetilde{k}_{1}-\widetilde{k}_{\overline{1}}.$

\noindent$\mathrm{(c):}$ Consider\textrm{ }$\gamma$ appearing in $S_{2}$ and
suppose that $\lambda_{\overline{i}}\leq\lambda_{\overline{i+1}}%
-k_{\overline{i+1}}$ for all $i=1,...,n-2$ and $\lambda_{\overline{1}%
}-k_{\overline{1}}\leq0$. Set $\widetilde{\gamma}=s_{0}\circ\gamma.$ Then
$\gamma_{\overline{s}}=\widetilde{\gamma}_{\overline{s}}$ for $s\neq1$ and
$\widetilde{\gamma}_{\overline{1}}=-\lambda_{\overline{1}}-k_{1}%
+k_{\overline{1}}-1$. Write $\widetilde{k}_{i}=k_{i},$ $\widetilde
{k}_{\overline{i}}=k_{\overline{i}}$ for all $i=2,...,n-1$ and set
$\widetilde{k}_{1}=k_{\overline{1}}-\lambda_{\overline{1}},$ $\widetilde
{k}_{\overline{1}}=k_{1}+\lambda_{\overline{1}}+1.$ We have $\widetilde{k}%
_{1}\geq0$, $\lambda_{\overline{i}}\leq\lambda_{\overline{i+1}}-\widetilde
{k}_{\overline{i+1}}$ for all $i=1,...,n-2$ and $\lambda_{\overline{1}%
}-\widetilde{k}_{\overline{1}}<0.\;$Moreover $\widetilde{k}_{\overline{1}%
}+\cdot\cdot\cdot+\widetilde{k}_{\overline{n-1}}+\widetilde{k}_{1}+\cdot
\cdot\cdot+\widetilde{k}_{n-1}=r$ (thus $\gamma$ appears in $S_{1})$ and
$\widetilde{\gamma}_{\overline{1}}=\lambda_{\overline{1}}+\widetilde{k}%
_{1}-\widetilde{k}_{\overline{1}}.$

\noindent$\mathrm{(d):}$\textrm{ }Now consider\textrm{ }$\gamma$ appearing in
$S_{1}$ or $S_{2}$ and suppose that $\lambda_{\overline{s}}\leq\lambda
_{\overline{s+1}}-k_{\overline{s+1}}$ for any $s\in\{1,...,n-2\}$,
$\lambda_{\overline{1}}-k_{\overline{1}}\geq0$ (resp. $\lambda_{\overline{1}%
}-k_{\overline{1}}>0)$ if $\gamma$ appears in $S_{1}$ (resp.\ in $S_{2})$ and
there exists $i\in\{1,...,n-2\}$ such that $\lambda_{\overline{i+1}%
}-k_{\overline{i+1}}<\lambda_{\overline{i}}+k_{i}-k_{\overline{i}}.$ Define
$\widetilde{\gamma}=s_{i}\circ\gamma$ as above.\ Set $\widetilde{k}%
_{\overline{i+1}}=k_{\overline{i+1}}$, $\widetilde{k}_{\overline{i}%
}=k_{\overline{i}}$, $\widetilde{k}_{i+1}=\lambda_{\overline{i}}%
-\lambda_{\overline{i+1}}-k_{\overline{i}}+k_{i}+k_{\overline{i+1}}-1$ and
$\widetilde{k}_{i}=(\lambda_{\overline{i+1}}-\lambda_{\overline{i}%
}-k_{\overline{i+1}})+k_{i+1}+k_{\overline{i}}+1.$ Write $\widetilde{k}%
_{t}=k_{t}$ for any $t\neq i,i+1,\overline{i},\overline{i+1}.$ We obtain
$\widetilde{k}_{i}\geq0$ and $\widetilde{k}_{i+1}\geq0$ since $\lambda
_{\overline{i}}\leq\lambda_{\overline{i+1}}-k_{\overline{i+1}}$ and
$\lambda_{\overline{i+1}}-k_{\overline{i+1}}<\lambda_{\overline{i}}%
+k_{i}-k_{\overline{i}}.$ Since $\widetilde{k}_{\overline{s}}=k_{\overline{s}%
}$ for all $s=1,...,n-1,$ we have $\lambda_{\overline{s}}\leq\lambda
_{\overline{s+1}}-\widetilde{k}_{\overline{s+1}}$ for any $s\in\{1,...,n-2\}$
and $\lambda_{\overline{1}}-\widetilde{k}_{\overline{1}}\geq0.\;$Moreover the
assertion $\lambda_{\overline{i+1}}-\widetilde{k}_{\overline{i+1}}%
<\lambda_{\overline{i}}+\widetilde{k}_{i}-\widetilde{k}_{\overline{i}}$ holds
since it is equivalent to $0<k_{i+1}+1.$ Finally $\widetilde{k}_{\overline{1}%
}+\cdot\cdot\cdot+\widetilde{k}_{\overline{n-1}}+\widetilde{k}_{1}+\cdot
\cdot\cdot+\widetilde{k}_{n-1}=k_{\overline{1}}+\cdot\cdot\cdot+k_{\overline
{n-1}}+k_{1}+\cdot\cdot\cdot+k_{n-1}$ and for any $s\in\{1,...,n-2\}$%
\[
\widetilde{\gamma}_{\overline{s}}=\lambda_{\overline{s}}+\widetilde{k}%
_{s}-\widetilde{k}_{\overline{s}}.
\]

\noindent Denote by $E_{a},E_{d}$ the sets of multi-indices $(k_{\overline{1}%
},...,k_{\overline{n-1}},k_{1},...,k_{n-1})$ such that $k_{\overline{1}}%
+\cdot\cdot\cdot+k_{\overline{n-1}}+k_{1}+\cdot\cdot\cdot+k_{n-1}=r$ and
satisfying respectively the assertions $\mathrm{(a)}$\textrm{,} $\mathrm{(d)}%
$\textrm{.\ }Let\textrm{ }$f$ be the map defined on $E_{a}\cup E_{d}$ by
\[
f(\gamma)=\widetilde{\gamma}.
\]
Then by the above arguments $f$ is a bijection which verifies $f(E_{a})=E_{a}$
and $f(E_{d})=E_{d}$.\ Now the pairing $\gamma\longleftrightarrow
\widetilde{\gamma}$ provides the cancellation of all the $s_{\gamma}$ with
$\gamma=(\lambda_{\overline{n-1}}+k_{n-1}-k_{\overline{n-1}},...,\lambda
_{\overline{1}}+k_{1}-k_{\overline{1}})$ such that $(k_{\overline{1}%
},...,k_{\overline{n}},k_{1},...,k_{n})\in E_{a}\cup E_{d}$ appearing in
$S_{1}.$ Indeed $s_{(\mu_{\overline{n}}+R,\gamma)}=-s_{(\mu_{\overline{n}%
}+R,\widetilde{\gamma})}.$ We obtain similarly the cancellation of all the
$s_{\gamma}$ such that $\gamma$ verifies the assertions $\mathrm{(a)}$ or
$\mathrm{(d)}$ appearing in $S_{2}.$

\noindent Now write $E_{b}$ (resp. $E_{c})$ for the set of multi-indices
$(k_{\overline{1}},...,k_{\overline{n-1}},k_{1},...,k_{n-1})$ such that
$k_{\overline{1}}+\cdot\cdot\cdot+k_{\overline{n-1}}+k_{1}+\cdot\cdot
\cdot+k_{n-1}=r$ (resp. $r-1)$ and satisfying assertion $\mathrm{(b)}$\textrm{
}(resp. $\mathrm{(c))}$\textrm{.\ }Let\textrm{ }$e$ be the map defined on
$E_{b}\cup E_{c}$ by
\[
e(\gamma)=\widetilde{\gamma}.
\]
Then $e$ is a bijection which verifies $e(E_{b})=E_{c}$ and $e(E_{c})=E_{b}$
and the $s_{\gamma}$ such that $\gamma$ verifies the assertions $\mathrm{(b)}$
or $\mathrm{(c)}$ cancel in $S_{R,r}.$ Finally by Corollary \ref{cor_pieriB}
and Remark $\mathrm{(iii)}$ following Corollary \ref{cor_pieriD} we obtain
\[
S_{R,r}=\sum_{\gamma\in P_{B_{n-1}}^{+},\lambda\in\left(  \gamma\otimes
r\right)  _{B_{n-1}}}b_{\gamma,r}^{\lambda}s_{(\mu_{\overline{n}}+R,\gamma).}%
\]
Note that this equality is also true when $R=r=0$ if we set $S_{0,0}%
=s_{(\mu_{\overline{n}},\lambda)}.$ Thus we have%
\[
Q_{\mu}=\sum_{R=0}^{+\infty}\sum_{\underset{r\equiv R\operatorname{mod}%
2}{0\leq r\leq R}}q^{R}\sum_{\gamma\in P_{B_{n-1}}^{+},\lambda\in\left(
\gamma\otimes r\right)  _{B_{n-1}}}b_{\gamma,r}^{\lambda}K_{\lambda
,\mu^{\prime}}^{B_{n-1}}(q)s_{(\mu_{\overline{n}}+R,\gamma)}%
\]
which is equivalent to (\ref{rec_mor_b}). So the theorem is proved.
\end{proof}

\begin{theorem}
\label{Th_rec_morrD}Let $\mu\in P_{D_{n}}^{+}.$ Then%
\begin{equation}
Q_{\mu}^{D_{n}}=\sum_{\gamma\in P_{D_{n-1}}^{+}}\sum_{R=0}^{+\infty}%
\sum_{r+2m=R}q^{R}\sum_{\lambda\in\left(  \gamma\otimes r\right)  _{n-1}%
}d_{\gamma,r}^{\lambda}K_{\lambda,\mu^{\prime}}^{D_{n-1}}(q)s_{(\mu
_{\overline{n}}+R,\gamma)}. \label{rec_mor_D}%
\end{equation}
\end{theorem}

\begin{proof}
Set $R_{\overline{i}}=R_{\varepsilon_{\overline{n}}-\varepsilon_{\overline{i}%
}}$ for $i=1,...,n-1$ and $R_{i}=R_{\varepsilon_{\overline{n}}+\varepsilon
_{\overline{i}}}$ for $i=1,...,n.$ We obtain as in proof of Theorem
\ref{Th_rec_morrB}%
\begin{multline*}
Q_{\mu}=\sum_{\lambda\in P_{D_{n-1}}^{+}}K_{\lambda,\mu^{\prime}}(q)\times\\
\left(  \sum_{R=0}^{+\infty}\sum_{\kappa_{\overline{1}}+k_{\overline{2}}%
+\cdot\cdot\cdot+k_{\overline{n-1}}+\kappa_{1}+k_{2}+\cdot\cdot\cdot
+k_{n-1}=R}q^{R}(R_{1})^{\kappa_{1}}(R_{\overline{1}})^{\kappa_{\overline{1}}%
}(R_{2})^{k_{2}}(R_{\overline{2}})^{k_{\overline{2}}}\cdot\cdot\cdot
(R_{n-1})^{k_{n-1}}(R_{\overline{n-1}})^{k_{\overline{n-1}}}s_{(\mu
_{\overline{n}},\lambda)}\right)  =\\
\sum_{R=0}^{+\infty}\sum_{\lambda\in P_{D_{n-1}}^{+}}q^{R}K_{\lambda
,\mu^{\prime}}(q)\sum_{\kappa_{\overline{1}}+k_{\overline{2}}+\cdot\cdot
\cdot+k_{\overline{n-1}}+\kappa_{1}+k_{2}+\cdot\cdot\cdot+k_{n-1}=R}%
s_{(\mu_{\overline{n}}+R,\lambda_{\overline{n-1}}+k_{n-1}-k_{\overline{n-1}%
},\cdot\cdot\cdot,\lambda_{\overline{2}}+k_{2}-k_{\overline{2}},\lambda
_{\overline{1}}+\kappa_{1}-\kappa_{\overline{1}})}.
\end{multline*}
Fix $\lambda,R$ and consider%
\[
S_{R}=\sum_{\kappa_{\overline{1}}+k_{\overline{2}}+\cdot\cdot\cdot
+k_{\overline{n-1}}+\kappa_{1}+k_{2}+\cdot\cdot\cdot+k_{n-1}=R}s_{(\mu
_{\overline{n}}+R,\lambda_{\overline{n-1}}+k_{n-1}-k_{\overline{n-1}}%
,\cdot\cdot\cdot,\lambda_{\overline{2}}+k_{2}-k_{\overline{2}},\lambda
_{\overline{1}}+\kappa_{1}-\kappa_{\overline{1}})}.
\]
Set $\gamma=(\lambda_{\overline{n-1}}+k_{n-1}-k_{\overline{n-1}},\cdot
\cdot\cdot,\lambda_{\overline{2}}+k_{2}-k_{\overline{2}},\lambda_{\overline
{1}}+\kappa_{1}-\kappa_{\overline{1}}).$

\noindent$\mathrm{(a):}$\textrm{ }Consider\textrm{ }$\gamma$ appearing in
$S_{R}$ and suppose that there exists $i\in\{2,...,n-2\}$ such that
$\lambda_{\overline{i}}>\lambda_{\overline{i+1}}-k_{\overline{i+1}}.$ Then we
associate a $\widetilde{\gamma}$ verifying $\lambda_{\overline{i}}%
>\lambda_{\overline{i+1}}-\widetilde{k}_{\overline{i+1}}$ to $\gamma$ as we
have done in case $\mathrm{(a)}$ of the above proof. This is possible since
$s_{i}=(\overline{i+1},\overline{i})(i,i+1)\in W_{D_{n}}.$

\noindent$\mathrm{(b):}$\textrm{ }Consider\textrm{ }$\gamma$ appearing in
$S_{R}$ such that $\lambda_{\overline{1}}>\lambda_{\overline{2}}%
-k_{\overline{2}}.$ We set $\widetilde{\gamma}=s_{1}\circ\gamma,$
$\widetilde{k}_{2}=\kappa_{1},$ $\widetilde{\kappa}_{1}=k_{2},$ $\widetilde
{k}_{\overline{2}}=\lambda_{\overline{2}}-\lambda_{\overline{1}}%
+\kappa_{\overline{1}}+1$ and $\widetilde{\kappa}_{\overline{1}}%
=\lambda_{\overline{1}}-\lambda_{\overline{2}}+k_{\overline{2}}-1.$ Then
$\gamma$ appears in $S_{R}$ and verifies $\lambda_{\overline{1}}%
>\lambda_{\overline{2}}-\widetilde{k}_{\overline{2}}$ whatever the sign of
$\lambda_{\overline{1}}.$

\noindent$\mathrm{(c):}$\textrm{ }Consider\textrm{ }$\gamma$ appearing in
$S_{R}$ such that $-\lambda_{\overline{1}}>\lambda_{\overline{2}}%
-k_{\overline{2}}.$ We set $\widetilde{\gamma}=s_{0}\circ\gamma,$
$\widetilde{k}_{2}=\kappa_{\overline{1}},$ $\widetilde{\kappa}_{\overline{1}%
}=k_{2},$ $\widetilde{k}_{\overline{2}}=\lambda_{\overline{2}}+\lambda
_{\overline{1}}+\kappa_{1}+1$ and $\widetilde{\kappa}_{1}=-\lambda
_{\overline{1}}-\lambda_{\overline{2}}+k_{\overline{2}}-1.$ Then $\gamma$
appears in $S_{R}$ and verifies $-\lambda_{\overline{1}}>\lambda_{\overline
{2}}-\widetilde{k}_{\overline{2}}$ whatever the sign of $\lambda_{\overline
{1}}.$

\noindent$\mathrm{(d):}$ Consider\textrm{ }$\gamma$ appearing in $S_{R}$ and
suppose that $\lambda_{\overline{s}}\leq\lambda_{\overline{s+1}}%
-k_{\overline{s+1}}$ for any $s\in\{2,...,n-2\}$, $\pm\lambda_{\overline{1}%
}>\lambda_{\overline{2}}-k_{\overline{2}},$ and there exists $i\in
\{1,...,n-2\}$ such that $\lambda_{\overline{i+1}}-k_{\overline{i+1}}%
<\lambda_{\overline{i}}+k_{i}-k_{\overline{i}}.$ We set $\widetilde{\gamma
}=s_{i}\circ\gamma$ and proceed as in case $\mathrm{(d)}$ of the above proof.

\noindent$\mathrm{(e):}$ Consider\textrm{ }$\gamma$ appearing in $S_{R}$ and
suppose that $\lambda_{\overline{s}}\leq\lambda_{\overline{s+1}}%
-k_{\overline{s+1}}$ for any $s\in\{2,...,n-2\}$, $\pm\lambda_{\overline{1}%
}>\lambda_{\overline{2}}-k_{\overline{2}},$ and $\lambda_{\overline{2}%
}-k_{\overline{2}}<\lambda_{\overline{1}}+\kappa_{1}-\kappa_{\overline{1}}.$
We set $\widetilde{\gamma}=s_{1}\circ\gamma$, $\widetilde{k}_{\overline{2}%
}=k_{\overline{2}},$ $\widetilde{\kappa}_{\overline{1}}=\kappa_{\overline{1}%
},$ $\widetilde{k}_{2}=\lambda_{\overline{1}}-\lambda_{\overline{2}}%
-\kappa_{\overline{1}}+\kappa_{1}-1$ and $\widetilde{\kappa}_{1}%
=\lambda_{\overline{2}}-\lambda_{\overline{1}}-k_{\overline{2}}+k_{2}%
+\kappa_{\overline{1}}+1.$ Then $\gamma$ appears in $S_{R}$ and verifies
$\lambda_{\overline{s}}\leq\lambda_{\overline{s+1}}-\widetilde{k}%
_{\overline{s+1}}$ for any $s\in\{2,...,n-2\}$, $\pm\lambda_{\overline{1}%
}>\lambda_{\overline{2}}-\widetilde{k}_{\overline{2}},$ and $\lambda
_{\overline{2}}-\widetilde{k}_{\overline{2}}<\lambda_{\overline{1}}%
+\widetilde{\kappa}_{1}-\widetilde{\kappa}_{\overline{1}}$ whatever the sign
of $\lambda_{\overline{1}}.$

\noindent$\mathrm{(f):}$ Consider\textrm{ }$\gamma$ appearing in $S_{R}$ and
suppose that $\lambda_{\overline{s}}\leq\lambda_{\overline{s+1}}%
-k_{\overline{s+1}}$ for any $s\in\{2,...,n-2\}$, $\pm\lambda_{\overline{1}%
}>\lambda_{\overline{2}}-k_{\overline{2}},$ and $\lambda_{\overline{2}%
}-k_{\overline{2}}<-\lambda_{\overline{1}}-\kappa_{1}+\kappa_{\overline{1}}.$
We set $\widetilde{\gamma}=s_{0}\circ\gamma$, $\widetilde{k}_{\overline{2}%
}=k_{\overline{2}},$ $\widetilde{\kappa}_{1}=\kappa_{1},$ $\widetilde{k}%
_{2}=-\lambda_{\overline{1}}-\lambda_{\overline{2}}+\kappa_{\overline{1}%
}-\kappa_{1}-1$ and $\widetilde{\kappa}_{1}=\lambda_{\overline{2}}%
+\lambda_{\overline{1}}-k_{\overline{2}}+k_{2}+\kappa_{1}+1.$ Then $\gamma$
appears in $S_{R}$ and verifies $\lambda_{\overline{s}}\leq\lambda
_{\overline{s+1}}-\widetilde{k}_{\overline{s+1}}$ for any $s\in\{2,...,n-2\}$,
$\pm\lambda_{\overline{1}}>\lambda_{\overline{2}}-\widetilde{k}_{\overline{2}%
},$ and $\lambda_{\overline{2}}-\widetilde{k}_{\overline{2}}<-\lambda
_{\overline{1}}+\widetilde{\kappa}_{1}-\widetilde{\kappa}_{\overline{1}}$
whatever the sign of $\lambda_{\overline{1}}.$

\noindent By considering the pairing $\gamma\longleftrightarrow\widetilde
{\gamma},$ the $s_{\gamma}$ appearing in $S_{R}$ cancel if they do not verify
simultaneously all the following conditions%
\begin{equation}
\left\{
\begin{tabular}
[c]{l}%
$1:\lambda_{\overline{1}}\leq\lambda_{\overline{2}}-\widetilde{k}%
_{\overline{2}}\text{ and }-\lambda_{\overline{1}}\leq\lambda_{\overline{2}%
}-\widetilde{k}_{\overline{2}}$\\
$2:\lambda_{\overline{i}}\leq\lambda_{\overline{i+1}}-k_{\overline{i+1}}$ for
$i=2,...,n-1$\\
$3:\lambda_{\overline{i+1}}-k_{\overline{i+1}}\geq\lambda_{\overline{i}}%
+k_{i}-k_{\overline{i}}$ for $i=2,...,n-2$\\
$4:\lambda_{\overline{2}}-k_{\overline{2}}\leq\lambda_{\overline{1}}%
+\kappa_{1}-\kappa_{\overline{1}}$ and $\lambda_{\overline{2}}-k_{\overline
{2}}\leq-\lambda_{\overline{1}}+\kappa_{1}-\kappa_{\overline{1}}$%
\end{tabular}
\right.  . \label{cond}%
\end{equation}
Note that conditions $1,2$ and $3$ are precisely conditions $\mathrm{(i)}%
,\mathrm{(ii)}$ and $\mathrm{(iii)}$ of Corollary \ref{cor_pieriD}. Let
$E_{R}$ be the set multi-indices $M=(\kappa_{\overline{1}},...,k_{\overline
{n}},\kappa_{1},...,k_{n})$ such that $\kappa_{\overline{1}}+\cdot\cdot
\cdot+k_{\overline{n-1}}+\kappa_{1}+\cdot\cdot\cdot+k_{n-1}=R$ and satisfying
(\ref{cond}). We can write $S=\sum_{M\in E_{R}}s_{(\mu_{\overline{n}}%
+R,\gamma_{M})}.$\textrm{ }Set $E_{R}^{-}=\{M\in E_{R},\kappa_{1}%
-\kappa_{\overline{1}}\leq0\}$ and $E_{R}^{+}=\{M\in E_{R},\kappa_{1}%
-\kappa_{\overline{1}}>0\}.$

\noindent Let $m$ be an integer such that $0\leq m\leq R/2.$ Set $r=R-2m.$
Consider the multi-indices $M\in E_{R}^{-}$ such that $\kappa_{1}=m$. Set
$k_{\overline{1}}=\kappa_{\overline{1}}-m=\kappa_{\overline{1}}-\kappa_{1}.$
If $\gamma_{\overline{1}}=\lambda_{\overline{1}}-k_{\overline{1}}\geq0$ (resp.
$\gamma_{\overline{1}}<0)$ then condition $4$ of (\ref{cond}) is equivalent to
condition $\mathrm{(iv,(a))}$ of Corollary \ref{cor_pieriD} (resp. to
condition $\mathrm{(iv,(d)).}$ Moreover $k_{\overline{1}}+\sum_{2\leq i\leq
n}(k_{\overline{i}}+k_{i})=r.$ Write $B^{-}(r)$ for the sub-graph of
$B^{D_{n-1}}(r)$ defined by the vertices which does not contain any letter
$1.$ Set $B(\gamma)\otimes B^{-}(r)=\underset{\lambda\in P_{D_{n-1}}^{+}}%
{\cup}B(\lambda)^{d_{\gamma,r}^{\lambda,-}}$ and $\left(  \gamma\otimes
r\right)  _{D_{n-1}}^{-}=\{\lambda\in P_{D_{n}}^{+},$ $d_{\gamma,r}%
^{\lambda,-}\neq0\}.$ By Remark $\mathrm{(iii)}$ following Corollary
\ref{cor_pieriD} we know that $\gamma\in P\in P_{D_{n-1}},$ so we obtain
\[
\sum_{M\in E_{R}^{-},\kappa_{1}=m}s_{(\mu_{\overline{n}}+R,\gamma_{M})}%
=\sum_{\gamma\in P_{D_{n-1},}\lambda\in\left(  \gamma\otimes r\right)
_{D_{n-1}}^{-}}d_{\gamma,r}^{\lambda,-}s_{(\mu_{\overline{n}}+R,\gamma)}.
\]
Now consider the multi-indices $M\in E_{R}^{+}$ such that $\kappa
_{\overline{1}}=m$. Set $r=R-2m$ and $B(\gamma)\otimes B^{+}(r)=\underset
{\lambda\in P_{D_{n-1}}^{+}}{\cup}B(\lambda)^{d_{\gamma,r}^{\lambda,+}}$ where
$B^{+}(r)$ is the sub-graph of $B^{D_{n-1}}(r)$ defined by the vertices which
does not contain any letter $\overline{1}.$ Write $\left(  \gamma\otimes
r\right)  _{D_{n-1}}^{+}=\{\lambda\in P_{D_{n}}^{+},$ $d_{\gamma,r}%
^{\lambda,+}\neq0\}.$ We obtain similarly%
\[
\sum_{M\in E_{R}^{+},\kappa_{\overline{1}}=m}s_{(\mu_{\overline{n}}%
+R,\gamma_{M})}=\sum_{\gamma\in P_{D_{n-1},}\lambda\in\left(  \gamma\otimes
r\right)  _{D_{n-1}}^{+}}d_{\gamma,r}^{\lambda,+}s_{(\mu_{\overline{n}%
}+R,\gamma)}.
\]
Finally%
\[
S=\sum_{r+2m=R}\ \sum_{\gamma\in P_{D_{n-1},},\lambda\in\left(  \gamma\otimes
r\right)  _{D_{n-1}}^{-}\cup\left(  \gamma\otimes r\right)  _{D_{n-1}}^{-}%
}\ (d_{\gamma,r}^{\lambda,-}+d_{\gamma,r}^{\lambda,+})s_{(\mu_{\overline{n}%
}+R,\gamma)}=\sum_{r+2m=R}\ \sum_{\gamma\in P_{D_{n-1},},\lambda\in\left(
\gamma\otimes r\right)  _{D_{n-1}}}d_{\gamma,r}^{\lambda}s_{(\mu_{\overline
{n}}+R,\gamma)}%
\]
since $\left(  \gamma\otimes r\right)  _{D_{n-1}}$ is the disjoint union of
$\left(  \gamma\otimes r\right)  _{D_{n-1}}^{-}$ and $\left(  \gamma\otimes
r\right)  _{D_{n-1}}^{+}.$ So the theorem is proved.
\end{proof}

\bigskip

Consider $\nu,\mu$ two generalized partitions of length $n.$ Write $p$ for the
lowest integer in $\{1,...,n\}$ such that $\nu_{\overline{p}}+p-\mu
_{\overline{n}}-n\geq0.$ For any $k\in\{p,p+1,...,n\}$ let $\sigma_{k}$ be the
signed permutation defined by%
\[
\sigma_{k}(i)=\left\{
\begin{tabular}
[c]{l}%
$i+1$ if $k\leq i\leq n-1$\\
$i$ if $1\leq i\leq k-1$\\
$k$ if $i=n$%
\end{tabular}
\right.  .
\]
Note that $(-1)^{l_{B}(\sigma_{k})}=(-1)^{l_{D}(\sigma_{k})}=(-1)^{n-k}.$ Let
$\gamma_{k}$ be the generalized partition of length $n-1$%
\[
\gamma_{k}=(\nu_{\overline{n}}+1,\nu_{\overline{n-1}}+1,...,\nu_{\overline
{k+1}}+1,\nu_{\overline{k-1}},...,\nu_{\overline{1}}).
\]
Finally set $R_{k}=\nu_{\overline{k}}+k-\mu_{\overline{n}}-n.$

\noindent From the above recurrence formulas it is possible to express any
Kostka-Foulkes polynomial $K_{\nu,\mu}(q)$ associated to a classical root
system of rank $n$ in terms of Kostka-Foulkes polynomials associated to the
corresponding root system of rank $n-1.$

\begin{theorem}
\label{Th_mor_expli}With the above notation we have
\begin{gather*}
\mathrm{(i)}:K_{\nu,\mu}^{B_{n}}(q)=\sum_{k=p}^{n}(-1)^{n-k}\times q^{R_{k}%
}\times\sum_{r+2m=R_{k}}\sum_{\lambda\in\left(  \gamma_{r}\otimes r\right)
_{B_{n-1}}}b_{\gamma_{r},r}^{\lambda}K_{\lambda,\mu^{\prime}}^{B_{n-1}}(q),\\
\mathrm{(ii)}:K_{\nu,\mu}^{C_{n}}(q)=\sum_{k=p}^{n}(-1)^{n-k}\times
\sum_{r+2m=R_{k}}\sum_{\lambda\in\left(  \gamma_{r}\otimes r\right)
_{C_{n-1}}}q^{R_{k}-m}\times c_{\gamma_{r},r}^{\lambda}K_{\lambda,\mu^{\prime
}}^{C_{n-1}}(q)\\
\mathrm{(iii)}:K_{\nu,\mu}^{D_{n}}(q)=\sum_{k=p}^{n}(-1)^{n-k}\times q^{R_{k}%
}\times\sum_{r+2m=R_{k}}\sum_{\lambda\in\left(  \gamma_{r}\otimes r\right)
_{D_{n-1}}}d_{\gamma_{r},r}^{\lambda}K_{\lambda,\mu^{\prime}}^{D_{n-1}}(q).
\end{gather*}
\end{theorem}

\begin{proof}
In case $\mathrm{(i),}$ write $E_{\nu}$ for the set of pairs $(\gamma,R)$ such
that there exists $\sigma_{(\gamma,R)}\in W_{B_{n}}$ verifying $\sigma
_{(\gamma,R)}\circ(\mu_{\overline{n}}+R,\gamma)=\nu.$ We obtain from Theorems
\ref{th_hall_kostka} and \ref{Th_rec_morrB}%
\begin{equation}
K_{\nu,\mu}(q)=\sum_{(\gamma,R)\in E_{\nu}}\sum_{r+2m=R}q^{R}\sum_{\lambda
\in\left(  \gamma\otimes r\right)  _{B_{n-1}}}b_{\gamma,r}^{\lambda
}(-1)^{l(\sigma_{(\gamma,R)})}K_{\lambda,\mu^{\prime}}(q).
\label{rec_mor_expli}%
\end{equation}
Consider $(\gamma,R)\in E_{\nu}.$ We must have%
\[
\sigma\left(  \mu_{\overline{n}}+R+n-\dfrac{1}{2},\gamma_{\overline{n-1}%
}+n-\dfrac{3}{2},...,\gamma_{\overline{1}}+\dfrac{1}{2}\right)  =\left(
\nu_{\overline{n}}+n-\dfrac{1}{2},\nu_{\overline{n-1}}+n-\dfrac{3}{2}%
,...,\nu_{\overline{1}}+\dfrac{1}{2}\right)  .
\]
The strictly decreasing subsequence $(\gamma_{\overline{n-1}}+n-\dfrac{3}%
{2},...,\gamma_{\overline{1}}+\dfrac{1}{2})$ must be sent under the action of
$\sigma$ on a strictly decreasing subsequence $I_{\gamma}$ of $(\nu
_{\overline{n}}+n-\dfrac{1}{2},\nu_{\overline{n-1}}+n-\dfrac{3}{2}%
,...,\nu_{\overline{1}}+\dfrac{1}{2}).$ These subsequences correspond to the
choice of a $\nu_{\overline{k}}+\dfrac{2k-1}{2}$ (for the image of
$\mu_{\overline{n}}+R+n-\dfrac{1}{2}$ under the action of $\sigma)$ which does
not belong to $I_{\gamma}.$ For such a subsequence we must have $\mu
_{\overline{n}}+R+n-\dfrac{1}{2}=\nu_{\overline{k}}+\dfrac{2k-1}{2}.$ Since
$R=\nu_{\overline{k}}+k-\mu_{\overline{n}}-n\geq0$ this implies that
$k\in\{p,...n\}$, $R=R_{k},$ $\sigma=\sigma_{k}$ and $\gamma=\gamma_{k}.$

\noindent We prove $\mathrm{(ii)}$ and $\mathrm{(iii)}$ similarly.
\end{proof}

\section{The statistics $\chi_{n}^{B},\chi_{n}^{C}$ and $\chi_{n}^{D}$}

In this section we introduce a statistic on partition shaped
Kashiwara-Nakashima's tableaux verifying%
\[
K_{\nu,\mu}(q)=\sum_{T\in\mathbf{T}(\lambda)_{\mu}}q^{\chi_{n}(T)}%
\]
when $(\nu,\mu)$ satisfies restrictive conditions. Although the statistic
$\chi_{n}$ can be regarded as a generalization of Lascoux-Sch\"{u}tzenberger's
charge for semi-standard tableaux, it does not permit to recover the
Kostka-Foulkes polynomial $K_{\nu,\mu}(q)$ for any $(\nu,\mu).$

\subsection{Catabolism}

\noindent From Theorem \ref{Th_mor_expli} we derive the following lemma:

\begin{lemma}
\label{lem_K(q)_perf}Let $\nu,\mu\in P_{n}^{+}$ be such that $\mu
_{\overline{n}}\geq\nu_{\overline{n-1}}.$ Set $l=\nu_{\overline{n}}%
-\mu_{\overline{n}}\geq0$ (otherwise $K_{\nu,\mu}(q)=0$).$\;$Then:%
\begin{align*}
\mathrm{(i)}:K_{\nu,\mu}^{B_{n}}(q)=q^{l}\sum_{r+2m=l}\text{\ }\sum
_{\lambda\in\left(  \nu^{\prime}\otimes r\right)  _{B_{n-1}}}b_{\nu^{\prime
},r}^{\lambda}K_{\lambda,\mu^{\prime}}^{B_{n-1}}(q),\\
\mathrm{(ii)}:K_{\nu,\mu}^{C_{n}}(q)=\sum_{r+2m=l}q^{r+m}\text{\ }%
\sum_{\lambda\in\left(  \nu^{\prime}\otimes r\right)  _{C_{n-1}}}%
c_{\nu^{\prime},r}^{\lambda}K_{\lambda,\mu^{\prime}}^{C_{n-1}}(q),\\
\mathrm{(iii)}:K_{\nu,\mu}^{D_{n}}(q)=q^{l}\sum_{r+2m=l}\text{\ }\sum
_{\lambda\in\left(  \nu^{\prime}\otimes r\right)  _{D_{n-1}}}d_{\nu^{\prime
},r}^{\lambda}K_{\lambda,\mu^{\prime}}^{D_{n-1}}(q).
\end{align*}
\end{lemma}

\begin{proof}
Assertions $\mathrm{(i),(ii)}$ and $\mathrm{(iii)}$ follow by applying Theorem
\ref{Th_mor_expli} with $p=n.$
\end{proof}

\noindent From now $\nu$\ and $\mu$ are generalized partitions with integers
parts. Consider $T\in\mathbf{T}^{n}(\nu)_{\mu}.$ Accordingly to Lemma
\ref{lem_fact_row}, we can write
\begin{equation}
\mathrm{w}(T)\equiv_{n}\mathrm{w}(R)\otimes\mathrm{w}(T^{\prime}).
\label{fact_tab}%
\end{equation}
Let $R^{\prime}$ be the row tableau obtained by erasing all the letters
$\overline{n}$ and $n$ in $R.$ The catabolism of the tableau $T$ is defined
by
\[
\mathrm{cat}(T)=P_{n-1}(\mathrm{w}(T^{\prime})\otimes\mathrm{w}(R^{\prime})).
\]
The tableau $\mathrm{cat}(T)$ is well defined and belongs to $\mathbf{T}%
^{n-1}(\lambda)_{\mu^{\prime}}$ where $\lambda$ is the shape of $\mathrm{cat}%
(T)$ since $T^{\prime}$ and $R^{\prime}$ do not contain any letter
$\overline{n}$ or $n.$

\noindent In the sequel we denote by $\mathrm{ch}_{A}$ the Lascoux-Sch\"{u}%
tzenberger's charge statistic on semi-standard tableaux. Note that
$\mathrm{ch}_{A}$ may be used to compute Kostka-Foulkes polynomials for the
root systems $B_{1}=C_{1}=A_{1}$ and $D_{3}=A_{3}.$

\noindent Consider $T\in\mathbf{T}^{n}(\nu)_{\mu}.$ The statistics $\chi
_{n}^{B},\chi_{n}^{C}$ and $\chi_{n}^{D}$ are defined recursively by:%
\begin{gather*}
\chi_{n}^{B}(T)=\left\{
\begin{tabular}
[c]{l}%
$\mathrm{ch}_{A}(T)$ if $n=1$\\
$\chi_{n-1}^{B}(\mathrm{cat}(T))+\nu_{\overline{n}}-\mu_{\overline{n}}$
otherwise
\end{tabular}
\right.  ,\text{ }\chi_{n}^{D}(T)=\left\{
\begin{tabular}
[c]{l}%
$\mathrm{ch}_{A}(T)$ if $n=3$\\
$\chi_{n-1}^{D}(\mathrm{cat}(T))+\nu_{\overline{n}}-\mu_{\overline{n}}$
otherwise
\end{tabular}
\right.  \text{ and}\\
\chi_{n}^{C}(T)=\left\{
\begin{tabular}
[c]{l}%
$\mathrm{ch}_{A}(T)$ if $n=1$\\
$\chi_{n-1}^{C}(\mathrm{cat}(T))+\nu_{\overline{n}}-\mu_{\overline{n}}-m$
otherwise
\end{tabular}
\right.  \text{ where }m\text{ is the number of letters }n\text{ in }R.
\end{gather*}

\noindent\textbf{Remark: }

\noindent$\mathrm{(i):}$ The statistics $\chi_{n}^{B},\chi_{n}^{C}$ and
$\chi_{n}^{D}$ can be regarded as extensions of $\mathrm{ch}_{A}$. More
precisely we have $\chi_{n}^{B}(T)=\chi_{n}^{C}(T)=\chi_{n}^{D}(T)=\mathrm{ch}%
_{A}(T)$ for the tableaux $T$ which contain only barred letters.

\noindent$\mathrm{(ii)}:$ To obtain $\chi_{1}^{B},\chi_{1}^{C}$ and $\chi
_{3}^{D}$ we need to compute $\mathrm{ch}_{A}$ on tableaux which are not
semi-standard. This can be done from the characterization of $\mathrm{ch}_{A}$
in terms of crystal graphs given in \cite{LLT} or more directly by using the
crystal graphs isomorphisms:%
\begin{equation}
B(\Lambda_{0}^{B_{1}})\simeq B(2\Lambda_{1}^{A_{1}})\text{, }B(\Lambda
_{0}^{C_{1}})\simeq B(\Lambda_{1}^{A_{1}}),\text{ }B(\Lambda_{0}^{D_{3}%
})\simeq B(\Lambda_{3}^{A_{3}})\text{, }B(\Lambda_{1}^{D_{3}})\simeq
B(\Lambda_{1}^{A_{3}})\text{ and }B(\Lambda_{2}^{D_{3}})\simeq B(\Lambda
_{2}^{A_{3}}) \label{isom}%
\end{equation}
which permit to turn each tableau $T$ related to types $B_{1},C_{1}$ and
$D_{3}$ into its corresponding tableau $\tau_{T}$ of type $A_{1}$ or $A_{3}$
via bumping algorithm on semi-standard tableaux.

\begin{example}
Consider the tableau of type $D_{3}$ and shape $(3,2,1),$ $T=%
\begin{tabular}
[c]{|l|ll}\hline
$\mathtt{\bar{3}}$ & $\mathtt{\bar{2}}$ & \multicolumn{1}{|l|}{$\mathtt{\bar
{1}}$}\\\hline
$\mathtt{1}$ & $\mathtt{2}$ & \multicolumn{1}{|l}{}\\\cline{1-1}\cline{1-2}%
$\mathtt{\bar{1}}$ &  & \\\cline{1-1}%
\end{tabular}
.$ Then $\mathrm{w}(T)=\overline{1}(\overline{2}2)(\overline{3}1\overline
{1}).$ We have $\overline{1}\in B(\Lambda_{2}^{D_{3}}),$ $(\overline{2}2)\in
B(\Lambda_{1}^{D_{3}}+\Lambda_{0}^{D_{3}})\simeq B(\Lambda_{1}^{A_{3}}%
+\Lambda_{3}^{A_{3}})$ and $(\overline{3}1\overline{1})\in B(2\Lambda
_{0}^{D_{3}})\simeq B(2\Lambda_{3}^{A_{3}}).$ Thus the semi-standard tableau
$\tau$ corresponding to $T$ is obtained by applying the bumping algorithm to
the word $w=(23)(3124)(124123).$ Finally $\tau_{T}=%
\begin{tabular}
[c]{|l|l|l|ll}\hline
$\mathtt{1}$ & $\mathtt{1}$ & $\mathtt{1}$ & $\mathtt{2}$ &
\multicolumn{1}{|l|}{$\mathtt{3}$}\\\hline
$\mathtt{2}$ & $\mathtt{2}$ & $\mathtt{2}$ & $\mathtt{3}$ &
\multicolumn{1}{|l}{}\\\cline{1-3}\cline{1-4}%
$\mathtt{3}$ & $\mathtt{4}$ & $\mathtt{4}$ &  & \\\cline{1-3}%
\end{tabular}
.$
\end{example}

\subsection{Catabolism and Kostka-Foulkes polynomials}

\noindent Consider $T\in\mathbf{T}^{n}(\nu)_{\mu}$ and suppose $n\geq2.$ For
any integer $p\leq n$ consider the sequence of tableaux defined by $T_{n}=T$
and $T_{k}=\mathrm{cat}(T_{k+1})$ for $k=n-1,...,p$. Denote by $v^{(k)}\in
P_{k}^{+}$ the shape of $T_{k}.$ Then $T_{k}\in\mathbf{T}^{k}(\nu^{(k)}%
)_{\mu^{(k)}}$ with $\mu^{(k)}=(\mu_{\overline{k}},...,\mu_{\overline{1}}).$

\begin{lemma}
\label{lem_tech}If $\mu_{\overline{p}}\geq v_{\overline{n-1}}$ then for every
$k=n,...,p$ we have $\mu_{\overline{p}}\geq v_{\overline{k-1}}^{(k)}.$
\end{lemma}

\begin{proof}
We proceed by induction on $k.$ The lemma is true for $k=n$. Consider
$k\in\{p+1,...,n\}$ such that $\mu_{\overline{p}}\geq v_{\overline{k-1}}%
^{(k)}.$ Then we must have $\nu_{\overline{k-2}}^{(k-1)}\leq\nu_{\overline
{k-1}}^{(k)}$ by Lemmas \ref{lem_b_gamm_tens_Lb} and \ref{lem_b_gamm_tens_Ld}
since the shape $\nu^{(k-1)}$ is obtained by adding or deleting boxes on
distinct columns of the shape obtained by deleting the longest row of
$\nu^{(k)}$. Hence $\nu_{\overline{k-2}}^{(k-1)}\leq\nu_{\overline{k-1}}%
^{(k)}\leq\mu_{\overline{p}}$.
\end{proof}

\bigskip

\begin{proposition}
\label{prop_xhi}Consider $\nu,\mu$ verifying one of the following conditions

\noindent$\mathrm{(i):}\nu,\mu\in P_{B_{n}}^{+}$ $n=1$ or, $n\geq2$ and
$\mu_{\overline{2}}\geq\nu_{\overline{n-1}}$

\noindent$\mathrm{(ii):}\nu,\mu\in P_{C_{n}}^{+}$ $n=1$ or, $n\geq2$ and
$\mu_{\overline{2}}\geq\nu_{\overline{n-1}}$

\noindent$\mathrm{(iii):}\nu,\mu\in P_{D_{n}}^{+}$ $n=3$ or, $n\geq4$ and
$\mu_{\overline{4}}\geq\nu_{\overline{n-1}}$

\noindent then%
\begin{equation}
K_{\nu,\mu}(q)=\sum_{T\in\mathbf{T}^{n}(\nu)_{\mu}}q^{\chi_{n}(T)}.
\label{K(q)_xhi}%
\end{equation}
\end{proposition}

\begin{proof}
The assertion is proved by induction on $n$.\ 

Case $\mathrm{(ii).}$ The proposition is true for the root system $C_{1}%
=A_{1}$. Now suppose that (\ref{K(q)_xhi}) is true for the root system
$C_{n-1}$ with $n\geq2$ and consider $\nu,\mu$ two partitions of length $n$
such that $\mu_{\overline{2}}\geq\nu_{\overline{n-1}}$. Set $l=\nu
_{\overline{n}}-\mu_{\overline{n}}.$ From Lemma \ref{lem_K(q)_perf}
$\mathrm{(i)}$ we obtain%
\[
K_{\nu,\mu}(q)=\sum_{r+2m=l}q^{r+m}\sum_{\lambda\in\left(  \nu^{\prime}\otimes
r\right)  _{C_{n-1}}}c_{\nu^{\prime},r}^{\lambda}K_{\lambda,\mu^{\prime}}(q)
\]
since $\mu_{\overline{n}}\geq\mu_{\overline{2}}\geq\nu_{\overline{n-1}}.$ Set
\[
K(q)=\sum_{T\in\mathbf{T}^{n}(\nu)_{\mu}}q^{\chi_{n}(T)}.
\]
Accordingly to Lemma \ref{lem_fact_row}, the reading of any $T\in
\mathbf{T}^{n}(\nu)_{\mu}$ can be factorized as
\[
\mathrm{w}(T)\equiv_{n}\mathrm{w}(R)\otimes\mathrm{w}(T^{\prime}).
\]
Set $\mathcal{T}_{m}=\{T\in\mathbf{T}^{n}(\nu)_{\mu},\mathrm{w}(R)$ contains
$m$ letters $n\}.$ We must have $0\leq m\leq l/2$ since all the letters
$\overline{n}$ or $n$ of $T$ belong to $R$ and the number of letters
$\overline{n}$ minus that of letters $n$ in $R$ must be equal to
$\mu_{\overline{n}}$. For any $T\in\mathcal{T}_{m}$ we can write
$\mathrm{cat}(T)=P_{n-1}(\mathrm{w}(T^{\prime})\otimes\mathrm{w}(R^{\prime}))$
where $R^{\prime}$ is a row tableau of length $r=l-2m$. The first row of
$T\ $contains at least $\mu_{\overline{n}}$ letters $\overline{n}$.\ Moreover
we have $\mu_{\overline{n}}\geq\mu_{\overline{2}}\geq\nu_{\overline{n-1}}$.
This means that $\{\mathrm{w}(R^{\prime})\otimes\mathrm{w}(T^{\prime}%
),T\in\mathcal{T}_{m}\}=\left(  B((r)_{n-1})\otimes B((\nu^{\prime})\right)
_{\mu^{\prime}}.$ Thus we have $\{\mathrm{w}(T^{\prime})\otimes\mathrm{w}%
(R^{\prime}),T\in\mathcal{T}_{m}\}=\left(  B(\nu^{\prime})\otimes
B((r)_{n-1})\right)  _{\mu^{\prime}}$ and $\{(\mathrm{cat}(T),T\in
\mathcal{T}_{m}\}$ is exactly the set of tableaux of shape $\lambda\in\left(
\nu^{\prime}\otimes r\right)  _{C_{n-1}}$ and weight $\mu^{\prime}.$ By lemma
\ref{lem_tech} we have $\mu_{\overline{2}}\geq\lambda_{\overline{n-2}}$ for
any $\lambda\in B(\nu^{\prime})\otimes B((r)_{n-1})$ when $n-1\geq2$.\ So we
can use the induction hypothesis and obtain
\begin{multline*}
K(q)=\sum_{m=0}^{l/2}\sum_{T\in\mathcal{T}_{m}}q^{\chi_{n}(T)}=\sum
_{r+2m=l}\sum_{T\in\mathcal{T}_{m}}q^{\chi_{n-1}(\mathrm{cat}(T))+l-m}%
=\sum_{r+2m=l}q^{r+m}\sum_{T\in\mathcal{T}_{m}}q^{\chi_{n-1}(\mathrm{cat}%
(T))}=\\
\sum_{r+2m=l}q^{r+m}\sum_{\lambda\in\left(  \nu^{\prime}\otimes r\right)
_{C_{n-1}}}c_{\nu^{\prime},r}^{\lambda}K_{\lambda,\mu^{\prime}}(q)=K_{\nu,\mu
}(q).
\end{multline*}

\noindent Assertions $\mathrm{(i)}$ and $\mathrm{(iii)}$ are proved similarly
by induction on $n$ starting respectively from $n=1$ and $n=3.$
\end{proof}

\begin{example}
Set $\nu=(4,1)$ and $\mu=(1,0)$ for type $B_{2}.$ For the $5$ corresponding
tableaux of shape $\lambda$ and weight $\mu$ we obtain:

$\chi_{2}^{B}\left(
\begin{tabular}
[c]{|l|lll}\hline
$\mathtt{\bar{2}}$ & $\mathtt{\bar{1}}$ & \multicolumn{1}{|l}{$\mathtt{\bar
{1}}$} & \multicolumn{1}{|l|}{$\mathtt{1}$}\\\hline
$\mathtt{1}$ &  &  & \\\cline{1-1}%
\end{tabular}
\right)  =\chi_{2}^{B}\left(  1\bar{1}\bar{1}\bar{2}\otimes1\right)
=\mathrm{ch}_{A}(1\otimes1\bar{1}\bar{1})+3=4+3=7,\vspace{0.1cm}$

$\chi_{2}^{B}\left(
\begin{tabular}
[c]{|l|lll}\hline
$\mathtt{\bar{2}}$ & $\mathtt{\bar{2}}$ & \multicolumn{1}{|l}{$\mathtt{0}$} &
\multicolumn{1}{|l|}{$\mathtt{2}$}\\\hline
$\mathtt{0}$ &  &  & \\\cline{1-1}%
\end{tabular}
\right)  =\chi_{2}^{B}\left(  20\bar{2}\bar{2}\otimes0\right)  =\mathrm{ch}%
_{A}(0\otimes0)+3=1+3=4,\vspace{0.1cm}$

$\chi_{2}^{B}\left(
\begin{tabular}
[c]{|l|lll}\hline
$\mathtt{\bar{2}}$ & $\mathtt{\bar{1}}$ & \multicolumn{1}{|l}{$\mathtt{0}$} &
\multicolumn{1}{|l|}{$\mathtt{1}$}\\\hline
$\mathtt{0}$ &  &  & \\\cline{1-1}%
\end{tabular}
\right)  =\chi_{2}^{B}\left(  10\bar{1}\bar{2}\otimes0\right)  =\mathrm{ch}%
_{A}(0\otimes10\bar{1})+3=3+3=6,\vspace{0.1cm}$

$\chi_{2}^{B}\left(
\begin{tabular}
[c]{|l|lll}\hline
$\mathtt{\bar{2}}$ & $\mathtt{\bar{2}}$ & \multicolumn{1}{|l}{$\mathtt{\bar
{1}}$} & \multicolumn{1}{|l|}{$\mathtt{2}$}\\\hline
$\mathtt{1}$ &  &  & \\\cline{1-1}%
\end{tabular}
\right)  =\chi_{2}^{B}\left(  2\bar{1}\bar{2}\bar{2}\otimes1\right)
=\mathrm{ch}_{A}(1\otimes\bar{1})+3=2+3=5$,$\vspace{0.1cm}$

$\chi_{2}^{B}\left(
\begin{tabular}
[c]{|l|lll}\hline
$\mathtt{\bar{2}}$ & $\mathtt{\bar{1}}$ & \multicolumn{1}{|l}{$\mathtt{1}$} &
\multicolumn{1}{|l|}{$\mathtt{1}$}\\\hline
$\mathtt{\bar{1}}$ &  &  & \\\cline{1-1}%
\end{tabular}
\right)  =\chi_{2}^{B}\left(  11\bar{1}\bar{2}\otimes\bar{1}\right)
=\mathrm{ch}_{A}(\bar{1}\otimes11\bar{1})+3=2+3=5\vspace{0.1cm}$

$\chi_{2}^{B}\left(
\begin{tabular}
[c]{|l|lll}\hline
$\mathtt{\bar{2}}$ & $\mathtt{\bar{2}}$ & \multicolumn{1}{|l}{$\mathtt{1}$} &
\multicolumn{1}{|l|}{$\mathtt{2}$}\\\hline
$\mathtt{\bar{1}}$ &  &  & \\\cline{1-1}%
\end{tabular}
\right)  =\chi_{2}^{B}\left(  21\bar{2}\bar{2}\otimes\bar{1}\right)
=\mathrm{ch}_{A}(\bar{1}\otimes1)+3=0+3=3.$

Finally $K_{\nu,\mu}^{B_{2}}(q)=q^{7}+q^{6}+2q^{5}+q^{4}+q^{3}.$
\end{example}

\noindent The following corollary makes clear $K_{\nu,\mu}(q)$ when $\nu$ is a
row partition.

\begin{corollary}
\label{cor_xhi_L}Let $\nu,\mu$ be two partitions such that $\nu$ is a row
partition and $\mu_{\overline{1}}\geq0$. Set $h_{n}(\mu)=\underset
{i=1}{\overset{n}{\sum}}(n-i)\mu\overline{_{i}}.$ Then for any $R\in
\mathbf{T}^{n}(\nu)_{\mu}$ we have

$\mathrm{(i):}$ $\chi_{n}^{B}(R)=h_{n}(\mu)+2\underset{i=1}{\overset{n}{\sum}%
}(n-i+1)k_{i}$ if $0\notin R$ and $\chi_{n}^{B}(R)=h_{n}(\mu)+2\underset
{i=1}{\overset{n}{\sum}}(n-i+1)k_{i}+n$ otherwise,

$\mathrm{(ii):}$ $\chi_{n}^{C}(R)=h_{n}(\mu)+\underset{i=1}{\overset{n}{\sum}%
}(2(n-i)+1)k_{i}$

$\mathrm{(iii):}$ $\chi_{n}^{D}(R)=h_{n}(\mu)+2\underset{i=2}{\overset{n}%
{\sum}}(n-i+1)k_{i}$

\noindent where $k_{i}$ is the number of letters $i$ which belong to $R.$
\end{corollary}

\begin{proof}
We proceed by recurrence on $n.$

\noindent Suppose first $n=1$ for cases $\mathrm{(i)}$ and $\mathrm{(ii)}$. We
deduce from proposition \ref{prop_degreeK} that $K_{\nu,\mu}^{C_{1}%
}(q)=q^{\tfrac{\nu-\mu}{2}}$ and $K_{\nu,\mu}^{B_{1}}(q)=q^{\nu-\mu}.$ Thus
$\chi_{1}^{C}(R)=\tfrac{\nu-\mu}{2}=k_{1}$,
\[
\chi_{1}^{B}(R)=\nu-\mu=\left\{
\begin{tabular}
[c]{l}%
$2k_{1}$ if $0\notin R$\\
$2k_{1}+1$ otherwise
\end{tabular}
\right.
\]
and the Corollary holds for $n=1.$ The rest of the proof is similar to that of
proposition 3.2.3 in \cite{lec3}.

\noindent Now suppose $n=3$ for case $\mathrm{(iii)}.$ We can write
\[
R=%
\begin{tabular}
[c]{|l|l|l|l|l|}\hline
$\overline{3}^{k_{\overline{3}}}$ & $\overline{2}^{k_{\overline{2}}}$ &
$\overline{1}^{k_{\overline{1}}}$ & $2^{k_{2}}$ & $3^{k_{3}}$\\\hline
\end{tabular}
\]
where
\begin{tabular}
[c]{|l|}\hline
$a^{k}$\\\hline
\end{tabular}
means that there are $k$ boxes containing the letter $a$ in $R.$ Then the
semi-standard tableau associated to $R$ by (\ref{isom}) is
\[
R_{A}=%
\begin{tabular}
[c]{|l|l|l|l|l|}\hline
$1^{k_{\overline{3}}}$ & $1^{k_{\overline{2}}}$ & $2^{k_{\overline{1}}}$ &
$2^{k_{2}}$ & $3^{k_{3}}$\\\hline
$2^{k_{\overline{3}}}$ & $3^{k_{\overline{2}}}$ & $3^{k_{\overline{1}}}$ &
$4^{k_{2}}$ & $4^{k_{3}}$\\\hline
\end{tabular}
.
\]
By using the definition of the charge for semi-standard tableaux one verifies
that $\mathrm{ch}(R_{A})=\mu_{\overline{2}}+2\mu_{\overline{1}}+2k_{3}%
+4k_{2}=\chi_{3}^{D}(R).$ Thus the corollary holds for $n=3$ and we terminate
as in proof of proposition 3.2.3 in \cite{lec3}.
\end{proof}

\bigskip

\noindent\textbf{Remarks: }

\noindent$\mathrm{(i):}$ Write $(r)$ for the row partition whose non zero part
is equal to $r.$ From Proposition \ref{prop_xhi} and Corollary \ref{cor_xhi_L}%
, we deduce that for any partition $\mu\in P_{+}$ we have $K_{(r),\mu
}(q)=q^{h_{n}(\mu)}\times K_{(l),0}(q)$ with $l=r-\left|  \mu\right|  .$ If
$l$ is even we obtain $K_{(l),0}^{B_{n}}(q)=q^{l/2}K_{(l),0}^{C_{n}}(q)$ since
the row tableaux of types $B_{n}$ and $C_{n}$ are then identical. Moreover the
map $t$ defined from $\mathbf{T}^{B_{n-1}}((l))$ to $\mathbf{T}^{D_{n}}((l))$
by changing each barred letter $\overline{x}$ (resp. unbarred letter $x$) of
$R$ into $\overline{x+1}$ (resp. $x+1)$ is a bijection. Hence we have
\[
K_{(l),0}^{D_{n}}(q)=\sum_{R\in\mathbf{T}^{D_{n}}((l))_{0}}q^{2\underset
{i=2}{\overset{n}{\sum}}(n-i+1)k_{i}}=\sum_{t^{-1}(R)\in\mathbf{T}^{B_{n-1}%
}((l))_{0}}q^{l+\underset{j=1}{\overset{n-1}{\sum}}2(n-1-j)k_{j}}%
=K_{(l),0}^{B_{n-1}}(q)=q^{l/2}K_{(l),0}^{C_{n-1}}(q).
\]

\noindent$\mathrm{(ii):}$ The statistic $\chi_{n}$ can not be used to compute
any Kostka-Foulkes polynomial. For type $C_{2},$ $\lambda=(3,1)$ and
$\mu=(0,0)$ we have $K_{\lambda,\mu}(q)=q^{5}+q^{4}+q^{3}.$ By considering the
$3$ tableaux of type $C_{2},$ shape $\lambda$ and weight $\mu$ we obtain%
\[
\chi_{2}^{C}\left(
\begin{tabular}
[c]{|l|ll}\hline
$\mathtt{\bar{1}}$ & $\mathtt{\bar{1}}$ & \multicolumn{1}{|l|}{1}\\\hline
$\mathtt{1}$ &  & \\\cline{1-1}%
\end{tabular}
\right)  =5,\text{ }\chi_{2}^{C}\left(
\begin{tabular}
[c]{|l|ll}\hline
$\mathtt{\bar{2}}$ & $\mathtt{\bar{1}}$ & \multicolumn{1}{|l|}{2}\\\hline
$\mathtt{1}$ &  & \\\cline{1-1}%
\end{tabular}
\right)  =3\text{ and }\chi_{2}^{C}\left(
\begin{tabular}
[c]{|l|ll}\hline
$\mathtt{\bar{2}}$ & $\mathtt{1}$ & \multicolumn{1}{|l|}{2}\\\hline
$\mathtt{\bar{1}}$ &  & \\\cline{1-1}%
\end{tabular}
\right)  =2
\]
and $K_{\lambda,\mu}(q)\neq q^{5}+q^{3}+q^{2}.$

\subsection{Cyclage graphs for the orthogonal root systems}

In \cite{lec3} we have introduced a (co)-cyclage graph structure on tableaux
of type $C.$ We are going to see that such a structure also exists for the
partition shaped tableaux of types $B$ and $D$. For any $n\geq1$ we embed the
finite alphabets $\mathcal{A}_{n}^{B},\mathcal{A}_{n}^{C}$ and $\mathcal{A}%
_{n}^{D}$ respectively into the infinite alphabets%
\begin{gather*}
\mathcal{A}_{\infty}^{B}=\{\cdot\cdot\cdot<\overline{n}<\cdot\cdot
\cdot<\overline{1}<0<1<\cdot\cdot\cdot<n<\cdot\cdot\cdot\}\\
\mathcal{A}_{\infty}^{C}=\{\cdot\cdot\cdot<\overline{n}<\cdot\cdot
\cdot<\overline{1}<1<\cdot\cdot\cdot<n<\cdot\cdot\cdot\}\\
\mathcal{A}_{\infty}^{D}=\{\cdot\cdot\cdot<\overline{n}<\cdot\cdot
\cdot<\overline{2}<%
\begin{tabular}
[c]{l}%
$\overline{1}$\\
$1$%
\end{tabular}
<2<\cdot\cdot\cdot<n<\cdot\cdot\cdot\}.
\end{gather*}
The vertices of the crystal $G_{\infty}^{B}=\underset{n\geq0}{%
{\textstyle\bigoplus}
}G_{n}^{B},G_{\infty}^{C}=\underset{n\geq0}{%
{\textstyle\bigoplus}
}G_{n}^{C}$ and $G_{\infty}^{D}=\underset{n\geq0}{%
{\textstyle\bigoplus}
}G_{n}^{D}$ can be regarded as the words respectively on $\mathcal{A}_{\infty
}^{B},\mathcal{A}_{\infty}^{C}$ and $\mathcal{A}_{\infty}^{D}$.$\;$The
congruences obtained by identifying the vertices of $G_{\infty}^{B},G_{\infty
}^{C}$ and $G_{\infty}^{D}$ equal up to the plactic relations of length $3$
are respectively denoted by $\equiv_{B},\equiv_{C}$ and $\equiv_{D}.$ Set
$\mathbf{T}^{B}=\underset{n\geq0}{\cup}\mathbf{T}_{n}^{B},$ $\mathbf{T}%
^{C}=\underset{n\geq0}{\cup}\mathbf{T}_{n}^{C}$ and $\mathbf{T}^{D}%
=\underset{n\geq0}{\cup}\mathbf{T}_{n}^{D}.$

\noindent By Remark $\mathrm{(iii)}$ before Lemma \ref{lem_fact_row}, there
exits a unique tableau $P(w)$ such that $w\equiv\mathrm{w}(P(w))$ computed
from $w$ without using contraction relation.

\noindent In the sequel $\mu$ is a partition with $n$ integers parts. A
tableau $T\in\mathbf{T}$ is of weight $\mathrm{wt}(T)=\mu$ if $T\in
\mathbf{T}_{m}$ with $m\geq n$, $d_{\overline{i}}=\mu_{\overline{i}}$ for
$1\leq i\leq n$ and $d_{\overline{i}}=0$ for $i>m.$ Set $\mathbf{T}^{B}%
[\mu]=\{T\in\mathbf{T}^{B}$ of weight $\mu\},$ $\mathbf{T}^{C}[\mu
]=\{T\in\mathbf{T}^{C}$ of weight $\mu\}$ and $\mathbf{T}^{D}[\mu
]=\{T\in\mathbf{T}^{D}$ of weight $\mu\}.$

\noindent Consider $T=C_{1}\cdot\cdot\cdot C_{r}\in\mathbf{T}_{\mu}$ with
$r>1$ columns.\ The cocyclage operation is authorized for $T$ if $T$ contains
at least a column with a letter $n$ or without letter $\overline{n}.$ In this
case, let $x$ be the rightmost letter of the longest row of $T.$ We can write
$\mathrm{w}(T)=x\mathrm{w}(T_{\ast})$ where $T_{\ast}\in\mathbf{T}$. Then we
set%
\[
U(T)=P(\mathrm{w}(T_{\ast})x).
\]
This means that $U(T)$ is obtained by column inserting $x$ in $T_{\ast}$
without using contraction relation.

\noindent\textbf{Remarks:}

\noindent$\mathrm{(i)}\mathbf{:}$ If $\mathrm{wt}(T)=0$ then the cocyclage
operation is always authorized.

\noindent$\mathrm{(ii)}\mathbf{:}$ By convention there is no cocyclage
operation on the columns.

\bigskip

\noindent We endow the set $\mathbf{T}[\mu]$ with a structure of graph by
drawing an array $T\rightarrow T^{\prime}$ if and only if the cocyclage
operation is authorized on $T$ and $U(T)=T^{\prime}.$ Write $\Gamma(T)$ for
the connected component containing $T.$

\begin{example}
\label{cont_ex_ch_B}For $\mu=(0,0,0)$ the following graphs are connected
components of $\mathbf{T}^{B}[\mu]:$%
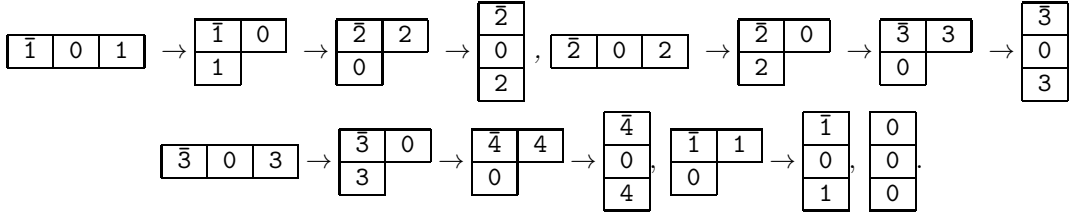
\begin{gather*}
\text{%
\begin{tabular}
[c]{|l|l|l|}\hline
$\mathtt{\bar{1}}$ & $\mathtt{0}$ & $\mathtt{1}$\\\hline
\end{tabular}
}\rightarrow\text{%
\begin{tabular}
[c]{|l|l}\hline
$\mathtt{\bar{1}}$ & \multicolumn{1}{|l|}{$\mathtt{0}$}\\\hline
$\mathtt{1}$ & \\\cline{1-1}%
\end{tabular}
}\rightarrow\text{%
\begin{tabular}
[c]{|l|l}\hline
$\mathtt{\bar{2}}$ & \multicolumn{1}{|l|}{$\mathtt{2}$}\\\hline
$\mathtt{0}$ & \\\cline{1-1}%
\end{tabular}
}\rightarrow\text{%
\begin{tabular}
[c]{|l|}\hline
$\mathtt{\bar{2}}$\\\hline
$\mathtt{0}$\\\hline
$\mathtt{2}$\\\hline
\end{tabular}
,
\begin{tabular}
[c]{|l|l|l|}\hline
$\mathtt{\bar{2}}$ & $\mathtt{0}$ & $\mathtt{2}$\\\hline
\end{tabular}
}\rightarrow\text{%
\begin{tabular}
[c]{|l|l}\hline
$\mathtt{\bar{2}}$ & \multicolumn{1}{|l|}{$\mathtt{0}$}\\\hline
$\mathtt{2}$ & \\\cline{1-1}%
\end{tabular}
}\rightarrow\text{%
\begin{tabular}
[c]{|l|l}\hline
$\mathtt{\bar{3}}$ & \multicolumn{1}{|l|}{$\mathtt{3}$}\\\hline
$\mathtt{0}$ & \\\cline{1-1}%
\end{tabular}
}\rightarrow\text{%
\begin{tabular}
[c]{|l|}\hline
$\mathtt{\bar{3}}$\\\hline
$\mathtt{0}$\\\hline
$\mathtt{3}$\\\hline
\end{tabular}
}\\%
\begin{tabular}
[c]{|l|l|l|}\hline
$\mathtt{\bar{3}}$ & $\mathtt{0}$ & $\mathtt{3}$\\\hline
\end{tabular}
\rightarrow%
\begin{tabular}
[c]{|l|l}\hline
$\mathtt{\bar{3}}$ & \multicolumn{1}{|l|}{$\mathtt{0}$}\\\hline
$\mathtt{3}$ & \\\cline{1-1}%
\end{tabular}
\rightarrow%
\begin{tabular}
[c]{|l|l}\hline
$\mathtt{\bar{4}}$ & \multicolumn{1}{|l|}{$\mathtt{4}$}\\\hline
$\mathtt{0}$ & \\\cline{1-1}%
\end{tabular}
\rightarrow%
\begin{tabular}
[c]{|l|}\hline
$\mathtt{\bar{4}}$\\\hline
$\mathtt{0}$\\\hline
$\mathtt{4}$\\\hline
\end{tabular}
,\text{ }%
\begin{tabular}
[c]{|l|l}\hline
$\mathtt{\bar{1}}$ & \multicolumn{1}{|l|}{$\mathtt{1}$}\\\hline
$\mathtt{0}$ & \\\cline{1-1}%
\end{tabular}
\rightarrow%
\begin{tabular}
[c]{|l|}\hline
$\mathtt{\bar{1}}$\\\hline
$\mathtt{0}$\\\hline
$\mathtt{1}$\\\hline
\end{tabular}
,\text{ }%
\begin{tabular}
[c]{|l|}\hline
$\mathtt{0}$\\\hline
$\mathtt{0}$\\\hline
$\mathtt{0}$\\\hline
\end{tabular}
.
\end{gather*}
All these tableaux belong to $T_{3}^{B}$ except
\begin{tabular}
[c]{|l|}\hline
$\mathtt{\bar{3}}$\\\hline
$\mathtt{0}$\\\hline
$\mathtt{3}$\\\hline
\end{tabular}
,
\begin{tabular}
[c]{|l|l}\hline
$\mathtt{\bar{3}}$ & \multicolumn{1}{|l|}{$\mathtt{0}$}\\\hline
$\mathtt{3}$ & \\\cline{1-1}%
\end{tabular}
,
\begin{tabular}
[c]{|l|l}\hline
$\mathtt{\bar{4}}$ & \multicolumn{1}{|l|}{$\mathtt{4}$}\\\hline
$\mathtt{0}$ & \\\cline{1-1}%
\end{tabular}
which belong to $T_{4}^{B}$ and
\begin{tabular}
[c]{|l|}\hline
$\mathtt{\bar{4}}$\\\hline
$\mathtt{0}$\\\hline
$\mathtt{4}$\\\hline
\end{tabular}
which belongs to $T_{5}^{B}.$
\end{example}

\noindent The following proposition is proved in the same way than Proposition
4.2.2 of \cite{lec3}.

\begin{proposition}
\label{prop_cyc}Let $T_{0}\in\mathbf{T}[0]$ and let $T_{k+1}=U(T_{k})$. Then
the sequence $(T_{n})$ is finite without repetition and there exists an
integer $e$ such that $T_{e}$ is a column of weight $0.$
\end{proposition}

\noindent In \cite{lec3} we introduce another statistic $\mathrm{ch}_{C_{n}}$
on Kashiwara-Nakashima's tableaux of type $C_{n}$ based on cocyclage
operation$.$ From $T\in\mathbf{T}^{C}[\mu]$ we define a finite sequence of
tableaux $(T_{k})_{0\leq k\leq p}$ whose last tableau $T_{p}$ is a column of
weight $0$. When $\mu=0$ this sequence $(T_{k})_{0\leq k\leq p}$ is precisely
that given in Proposition \ref{prop_cyc}. Then the statistic $\mathrm{ch}%
_{C_{n}}$ is first defined on the columns of weight $0$ next on the tableaux
by setting
\[
\mathrm{ch}_{C_{n}}(T)=\mathrm{ch}_{C_{n}}(C_{T})+p.
\]
We conjecture that (\ref{K(q)_xhi}) holds if we replace $\chi_{n}^{C}$ by
$\mathrm{ch}_{C_{n}}$ whatever the partitions $\lambda$ and $\mu.$ In
particular $\mathrm{ch}_{C_{n}}(T)\neq\chi_{n}(T)$ in general.

\noindent Unfortunately such a statistic defined in the same way for computing
Kostka-Foulkes polynomials can not exist for the orthogonal root systems. This
can be verified by considering the case $\left|  \lambda\right|  =3,$ $\mu=0$
for type $B_{3}$. Set $\lambda_{1}=(3,0,0),$ $\lambda_{2}=(2,1,0)$ and
$\lambda_{3}=(1,1,1).$ We have $K_{\lambda_{1},0}^{B_{3}}(q)=q^{9}+q^{7}%
+q^{5},$ $K_{\lambda_{2},0}^{B_{3}}(q)=q^{8}+q^{7}+q^{6}+q^{5}+q^{4}$ and
$K_{\lambda_{3},0}^{B_{3}}(q)=q^{6}+q^{4}+q^{2}.$ Then it is impossible to
associate a statistic $\mathrm{ch}_{B_{n}}$ to the $11$ tableaux of type
$B_{3},$ weight $0$ and shape $\lambda_{1},\lambda_{2}$ or $\lambda_{3}$
compatible with the cyclage graph structure given in Example
\ref{cont_ex_ch_B} (that is, such that $\mathrm{ch}_{B_{n}}(T)=\mathrm{ch}%
_{B_{n}}(T^{\prime})+1$ if $T\rightarrow T^{\prime}$) and relevant for
computing the corresponding Kostka-Foulkes polynomials. The situation is
similar for type $D_{3},$ $\left|  \lambda\right|  =3$ and $\mu=(1,0,0).$

\section{Explicit formulas for $K_{\lambda,\mu}(q)$}

\subsection{Explicit formulas for $\left|  \lambda\right|  \leq3$}

\noindent In the sequel we suppose that $\lambda$ is a partition such that
$\lambda_{\overline{1}}\geq0.$ We give below the matrix $K(q)=(K_{\lambda,\mu
}(q))$ with $\left|  \lambda\right|  \leq3$ associated to each root system
$B_{n},C_{n}$ and $D_{n}.$ When $\left|  \lambda\right|  =\left|  \mu\right|
,$ $K_{\lambda,\mu}(q)$ can be regarded as a Kostka-Foulkes polynomial for the
root system $A_{n-1}$. Such polynomials have been already compute (see
\cite{mac} p 329). So we only give the entries of $K(q)$ corresponding to a
weight $\mu$ such that $\left|  \mu\right|  \leq2$. In the following matrices
we have labelled the columns by $\lambda$ and the rows by $\mu$ and represent
each partition by its Young diagram. The expressions for the Kostka-Foulkes
polynomials are obtained by using Proposition \ref{prop_degreeK}, Theorem
\ref{Th_mor_expli}, Proposition \ref{prop_xhi} and Corollary \ref{cor_xhi_L}.

\subsubsection{$K(q)$-matrix for the root system $B_{n}$}%

\[%
\begin{array}
[c]{ccccccc}%
\vspace{0.2cm} &
\begin{tabular}
[c]{|l|l|l|}\hline
&  & \\\hline
\end{tabular}
&
\begin{tabular}
[c]{|l|l}\hline
& \multicolumn{1}{|l|}{}\\\hline
& \\\cline{1-1}%
\end{tabular}
&
\begin{tabular}
[c]{|l|}\hline
\\\hline
\\\hline
\\\hline
\end{tabular}
&
\begin{tabular}
[c]{|l|l|}\hline
& \\\hline
\end{tabular}
&
\begin{tabular}
[c]{|l|}\hline
\\\hline
\\\hline
\end{tabular}
&
\begin{tabular}
[c]{|l|}\hline
\\\hline
\end{tabular}
\\%
\begin{tabular}
[c]{|l|l|}\hline
& \\\hline
\end{tabular}
\vspace{0.2cm} & q^{n} & q^{n-1} & 0 & 1 & 0 & 0\\%
\begin{tabular}
[c]{|l|}\hline
\\\hline
\\\hline
\end{tabular}
\vspace{0.2cm} & q^{n+1} & q^{n}+q^{n-1} & q^{n-2} & q & 1 & 0\\%
\begin{tabular}
[c]{|l|}\hline
\\\hline
\end{tabular}
\vspace{0.4cm} & q^{2}\times\tfrac{q^{2n}-1}{q^{2}-1} & q^{n}+q\times
\tfrac{q^{2n-1}-1}{q-1} & q\times\tfrac{q^{2n-2}-1}{q^{2}-1} & q^{n} & q^{n-1}%
& 1\\
\emptyset &  q^{n+2}\times\tfrac{q^{2n-1}-1}{q-1} & q^{n+1}\times
\tfrac{q^{2n-1}-1}{q-1} & q^{n-1}\times\tfrac{q^{2n}-1}{q^{2}-1} & q^{2}%
\times\tfrac{q^{2n}-1}{q^{2}-1} & q\times\tfrac{q^{2n}-1}{q^{2}-1} & q^{n}%
\end{array}
\]

\subsubsection{$K(q)$-matrix for the root system $C_{n}$}%

\[%
\begin{array}
[c]{cccccc}%
\vspace{0.2cm} &
\begin{tabular}
[c]{|l|l|l|}\hline
&  & \\\hline
\end{tabular}
&
\begin{tabular}
[c]{|l|l}\hline
& \multicolumn{1}{|l|}{}\\\hline
& \\\cline{1-1}%
\end{tabular}
&
\begin{tabular}
[c]{|l|}\hline
\\\hline
\\\hline
\\\hline
\end{tabular}
&
\begin{tabular}
[c]{|l|l|}\hline
& \\\hline
\end{tabular}
&
\begin{tabular}
[c]{|l|}\hline
\\\hline
\\\hline
\end{tabular}
\\%
\begin{tabular}
[c]{|l|l|}\hline
& \\\hline
\end{tabular}
\vspace{0.2cm} & 0 & 0 & 0 & 1 & 0\\%
\begin{tabular}
[c]{|l|}\hline
\\\hline
\\\hline
\end{tabular}
\vspace{0.2cm} & 0 & 0 & 0 & q & 1\\%
\begin{tabular}
[c]{|l|}\hline
\\\hline
\end{tabular}
\vspace{0.4cm} & q\times\tfrac{q^{2n}-1}{q^{2}-1} & q\times\tfrac{q^{2n-2}%
-1}{q-1} & q^{2}\times\tfrac{q^{2n-4}-1}{q^{2}-1} & 0 & 0\\
\emptyset & 0 & 0 & 0 & q\times\tfrac{q^{2n}-1}{q^{2}-1} & q^{2}\times
\tfrac{q^{2n-2}-1}{q^{2}-1}%
\end{array}
\]

\subsubsection{$K(q)$-matrix for the root system $D_{n}$}%

\[%
\begin{array}
[c]{cccccc}%
\vspace{0.2cm} &
\begin{tabular}
[c]{|l|l|l|}\hline
&  & \\\hline
\end{tabular}
&
\begin{tabular}
[c]{|l|l}\hline
& \multicolumn{1}{|l|}{}\\\hline
& \\\cline{1-1}%
\end{tabular}
&
\begin{tabular}
[c]{|l|}\hline
\\\hline
\\\hline
\\\hline
\end{tabular}
&
\begin{tabular}
[c]{|l|l|}\hline
& \\\hline
\end{tabular}
&
\begin{tabular}
[c]{|l|}\hline
\\\hline
\\\hline
\end{tabular}
\\%
\begin{tabular}
[c]{|l|l|}\hline
& \\\hline
\end{tabular}
\vspace{0.2cm} & 0 & 0 & 0 & 1 & 0\\%
\begin{tabular}
[c]{|l|}\hline
\\\hline
\\\hline
\end{tabular}
\vspace{0.2cm} & 0 & 0 & 0 & q & 1\\%
\begin{tabular}
[c]{|l|}\hline
\\\hline
\end{tabular}
\vspace{0.4cm} & q^{2}\times\tfrac{q^{2n-2}-1}{q^{2}-1} & q^{n-1}%
+q\times\tfrac{q^{2n-3}-1}{q-1} & q^{n-2}+q\times\tfrac{q^{2n-4}-1}{q^{2}-1} &
0 & 0\\
\emptyset & 0 & 0 & 0 & q^{2}\times\tfrac{q^{2n-2}-1}{q^{2}-1} &
q^{n-1}+q\times\tfrac{q^{2n-2}-1}{q^{2}-1}%
\end{array}
\]

\noindent\textbf{Remark: }For $n\geq4$ the partitions $\lambda$ and $\mu$ in
the above matrix verify $\lambda^{\ast}=\lambda$ and $\mu^{\ast}=\mu.$ Hence
by (\ref{K=K*}) we have $K_{\lambda,\mu}^{D_{n}}(q)=K_{\lambda^{\ast}%
,\mu^{\ast}}^{D_{n}}(q)=K_{\lambda^{\ast},\mu}^{D_{n}}(q)=K_{\lambda,\mu
^{\ast}}^{D_{n}}(q).$

\subsection{Explicit formulas for the root system $B_{2}=C_{2}$ and $\mu=0$}

\noindent Note first that the roots systems $B_{2}$ and $C_{2}$ are identical.
More precisely denote by $\Psi$ the linear map%
\[
\Psi:\left\{
\begin{tabular}
[c]{c}%
$P_{B_{2}}^{+}\rightarrow P_{C_{2}}^{+}$\\
$(\lambda_{\overline{2}},\lambda_{\overline{1}})\longmapsto(\lambda
_{\overline{2}}+\lambda_{\overline{1}},\lambda_{\overline{2}}-\lambda
_{\overline{1}})$%
\end{tabular}
\right.  .
\]
Accordingly to (\ref{simple_roots}), the simple roots for the roots systems
$B_{2}$ and $C_{2}$ are $\alpha_{0}^{B_{2}}=\varepsilon_{\overline{1}}%
,\alpha_{1}^{B_{2}}=\varepsilon_{\overline{2}}-\varepsilon_{\overline{1}}$ and
$\alpha_{0}^{C_{2}}=2\varepsilon_{\overline{1}},\alpha_{1}^{C_{2}}%
=\varepsilon_{\overline{2}}-\varepsilon_{\overline{1}}.$ Thus we have
$\Psi(\alpha_{0}^{B_{2}})=\alpha_{1}^{C_{2}}$ and $\Psi(\alpha_{1}^{B_{2}%
})=\alpha_{0}^{C_{2}}.$ This implies the equality%
\begin{equation}
K_{(\lambda,\mu)}^{B_{2}}(q)=K_{\Psi(\lambda,\mu)}^{C_{2}}(q). \label{B2=C2}%
\end{equation}
So it is sufficient to explicit the Kostka-Foulkes polynomials for the root
system $C_{2}.$

\begin{proposition}
Let $\lambda=(\lambda_{\overline{2}},\lambda_{\overline{1}})$ be a generalized
partition of length $2.$

\begin{enumerate}
\item  If $\lambda\in P_{+}^{C_{2}}$ then
\[
K_{\lambda,0}^{C_{2}}(q)=\left\{
\begin{tabular}
[c]{l}%
$q^{\tfrac{\lambda_{\overline{2}}+\lambda_{\overline{1}}}{2}}\left(
\dfrac{q^{\lambda_{\overline{1}}+2}-1}{q^{2}-1}+q^{2}\times\dfrac
{q^{\lambda_{\overline{1}}+1}-1}{q-1}\times\dfrac{q^{\lambda_{\overline{2}%
}-\lambda_{\overline{1}}}-1}{q^{2}-1}\right)  $ if $\lambda_{\overline{2}}$
and $\lambda_{\overline{1}}$ are even\vspace{0.1cm}\\
$q^{\tfrac{\lambda_{\overline{2}}+\lambda_{\overline{1}}}{2}+1}\left(
\dfrac{q^{\lambda_{\overline{1}}+1}-1}{q^{2}-1}+q\times\dfrac{q^{\lambda
_{\overline{1}}+1}-1}{q-1}\times\dfrac{q^{\lambda_{\overline{2}}%
-\lambda_{\overline{1}}}-1}{q^{2}-1}\right)  $ if $\lambda_{\overline{2}}$ and
$\lambda_{\overline{1}}$ are odd\vspace{0.1cm}\\
$0$ otherwise.
\end{tabular}
\right.  .
\]

\item  If $\lambda\in P_{+}^{B_{2}}$ then
\[
K_{\lambda,0}^{B_{2}}(q)=\left\{
\begin{tabular}
[c]{l}%
$q^{\lambda_{\overline{2}}}\left(  \dfrac{q^{2\lambda_{\overline{1}}+2}%
-1}{q^{2}-1}+q^{2}\times\dfrac{q^{2\lambda_{\overline{1}}+1}-1}{q-1}%
\times\dfrac{q^{\lambda_{\overline{2}}-\lambda_{\overline{1}}}-1}{q^{2}%
-1}\right)  $ if $\lambda_{\overline{2}}+\lambda_{\overline{1}}$ is
even\vspace{0.1cm}\\
$q^{\lambda_{\overline{2}}+1}\times\dfrac{q^{2\lambda_{\overline{1}}+1}%
-1}{q-1}\times\dfrac{q^{\lambda_{\overline{2}}-\lambda_{\overline{1}}+1}%
-1}{q^{2}-1}$ otherwise
\end{tabular}
\right.  .
\]
\end{enumerate}
\end{proposition}

\begin{proof}
$1:$ Note first that $K_{\lambda,0}^{C_{2}}(q)=0$ if $\left|  \lambda\right|
$ is odd since all the tableaux of weight $0$ and type $C_{2}$ must have a
pair number of boxes. So we can suppose that $\lambda_{\overline{2}}$ and
$\lambda_{\overline{1}}$ have the same parity. By Theorem \ref{Th_mor_expli}
we must have
\[
K_{\lambda,0}^{C_{2}}(q)=\sum_{r+2m=\lambda_{\overline{2}}}q^{r+m}\sum
_{\eta\in((\lambda_{\overline{1}})\otimes r)_{1}}c_{(\lambda_{\overline{1}%
}),r}^{\eta}K_{\eta,0}^{C_{1}}(q)-\sum_{r+2m=\lambda_{\overline{1}}-1}%
q^{r+m}\sum_{\eta\in((\lambda_{\overline{2}}+1)\otimes r)_{1}}c_{(\lambda
_{\overline{2}}+1),r}^{\eta}K_{\eta,0}^{C_{1}}(q)
\]
where by abuse of notation the second sum is equal to $0$ if $\lambda
_{\overline{1}}=0.$ Now the $K_{\eta,0}^{C_{1}}(q)$'s are Kostka-Foulkes
polynomials for the root system $C_{1}=A_{1}$ hence $K_{\eta,0}^{C_{1}%
}(q)=q^{\eta/2}.$ Moreover Lemma \ref{lem_plu_hp} implies that
\[
B(\gamma)\otimes B(r)=\underset{p=0}{\overset{\min(\gamma,r)}{\cup}}%
B(\gamma+r-2p)
\]
for any integers $\gamma,r.$ We obtain%
\[
K_{\lambda,0}^{C_{2}}(q)=\sum_{r+2m=\lambda_{\overline{2}}}\sum_{p=0}%
^{\min(\lambda_{\overline{1}},r)}q^{r+m}\times q^{\tfrac{\lambda_{\overline
{1}}+r}{2}-p}-\sum_{r+2m=\lambda_{\overline{1}}-1}\sum_{p=0}^{r}q^{r+m}\times
q^{\tfrac{\lambda_{\overline{2}}+r+1}{2}-p}.
\]
Indeed we have $\min(\lambda_{\overline{2}}-1,r)=r$ in the second sum since
$r\leq\lambda_{\overline{1}}-1<\lambda_{\overline{2}}+1.$ This can be
rewritten as%
\begin{multline*}
K_{\lambda,0}^{C_{2}}(q)=\sum_{\underset{r\equiv\lambda_{\overline{2}}\text{
}\operatorname{mod}2}{r=0}}^{\lambda_{\overline{1}}}\sum_{p=0}^{r}%
q^{r+\tfrac{\lambda_{\overline{2}}-r}{2}+\tfrac{\lambda_{\overline{1}}+r}%
{2}-p}+\sum_{\underset{r\equiv\lambda_{\overline{2}}\text{ }\operatorname{mod}%
2}{r=\lambda_{\overline{1}}+1}}^{\lambda_{\overline{2}}}\sum_{p=0}%
^{\lambda_{\overline{1}}}q^{r+\tfrac{\lambda_{\overline{2}}-r}{2}%
+\tfrac{\lambda_{\overline{1}}+r}{2}-p}-\sum_{\underset{r\equiv\lambda
_{\overline{1}}-1\text{ }\operatorname{mod}2}{r=0}}^{\lambda_{\overline{1}}%
-1}\sum_{p=0}^{r}q^{r+\tfrac{\lambda_{\overline{1}}-r-1}{2}+\tfrac
{\lambda_{\overline{2}}+r+1}{2}-p}\\
=q^{\tfrac{\lambda_{\overline{2}}+\lambda_{\overline{1}}}{2}}\left(
\sum_{\underset{r\equiv\lambda_{\overline{2}}\text{ }\operatorname{mod}2}%
{r=0}}^{\lambda_{\overline{1}}}\sum_{p=0}^{r}q^{r-p}+\sum_{\underset
{r\equiv\lambda_{\overline{2}}\text{ }\operatorname{mod}2}{r=\lambda
_{\overline{1}}+1}}^{\lambda_{\overline{2}}}\sum_{p=0}^{\lambda_{\overline{1}%
}}q^{r-p}-\sum_{\underset{r\equiv\lambda_{\overline{1}}-1\text{ }%
\operatorname{mod}2}{r=0}}^{\lambda_{\overline{1}}-1}\sum_{p=0}^{r}%
q^{r-p}\right)  .
\end{multline*}
Then the Proposition easily follows by distinguishing the two cases
$\lambda_{\overline{2}}$ even and $\lambda_{\overline{2}}$ odd.

\noindent$2:$ This is an immediate consequence of $1$ and (\ref{B2=C2}).
\end{proof}

\noindent\textbf{Remark: }

\noindent$\mathrm{(i):}$ Similar formulas also exist for the root system
$A_{2}.$ For any partition $\lambda=(a,b,0)$ we have%
\[
K_{\lambda,0}^{A_{2}}(q)=\left\{
\begin{tabular}
[c]{l}%
$q^{a-b}\times\dfrac{q^{a+1}-1}{q-1}$ if $a\geq2b$\\
$q^{b}\times\dfrac{q^{a-b+1}-1}{q-1}$ otherwise.
\end{tabular}
\right.  .
\]

\bigskip

\noindent$\mathrm{(ii):}$ For a weight $\mu\neq0,$ the situation becomes more
complex and simple formulas for the $K_{\lambda,\mu}(q)$ seem do not exist.

\end{document}